\documentclass[aps,prx,reprint,superscriptaddress,floatfix]{revtex4-2}

\usepackage{amsmath,amssymb,amsfonts,bm,mathtools}
\usepackage{amsthm}
\usepackage{graphicx}
\usepackage{xcolor}
\usepackage{xr-hyper}
\usepackage{subcaption}
\usepackage[colorlinks=true,linkcolor=blue,citecolor=blue]{hyperref}
\DeclareMathOperator{\Tr}{Tr}
\newcommand{\R}{\mathbb{R}}
\newcommand{\E}{\mathbb{E}}
\newcommand{\Lf}{\mathcal{L}}
\newcommand{\Lb}{\mathcal{L}^{\dagger}}
\newcommand{\Kop}{\mathcal{K}}
\newcommand{\Atil}{\widetilde{A}}
\newcommand{\Stil}{\widetilde{\Sigma}}
\newcommand{\rhotil}{\widetilde{\rho}}
\newcommand{\Lp}{L_{+}}
\newcommand{\Lcal}{\Lambda_{+}}
\newcommand{\He}{\mathrm{He}}

\newtheorem{proposition}{Proposition}

\begin{document}

\title{\texorpdfstring{Local Koopman modes of saddle-point diffusions:\\ exact singular functions and a data-driven estimator}{title}}
\author{Xuehao Ding}
\affiliation{Flatiron Institute, Simons Foundation}

\begin{abstract}
Coherent sets—regions of phase space that resist mixing—organize transport in stochastic dynamical systems and are encoded in the leading singular functions of the stochastic Koopman operator. These are typically computed numerically and globally; a local, analytic theory near a fixed point has been missing, although fixed points and the transitions between them organize the flow. We supply the local theory for a saddle, the transition state through which probability flows between basins. Naively the construction fails—probability escapes along unstable directions, the invariant measure is non-normalizable, the singular functions leave $L^2$, and singular values exceed one. We show that conditioning on non-escape (a Doob's $h$-transform) repairs this; that the required $h$ is the unique stabilizing solution of an algebraic Riccati equation; and that the resulting slow modes decay at the escape rate along the unstable manifold. We obtain these modes in closed form—products of Hermite polynomials in canonical coordinates times the $h$-function—and use their unitary structure to derive a data-driven estimator: canonical correlation analysis on conditioned trajectories. When a single unstable direction dominates, the leading mode is the reaction coordinate, locally recovering the stable manifold and committor of transition-path theory. The construction connects coherent sets to transition-path theory, to the stabilizing Riccati solution of linear control, and to the time-reversed dynamics of generative diffusion models.
\end{abstract}

\maketitle

\section{Introduction}

A nonlinear dynamical system is organized by its fixed points and by the transitions between them: trajectories dwell near stable equilibria and cross between them through the neighborhoods of unstable equilibria~\cite{strogatz2024nonlinear, freidlin2012random}. Local behavior around fixed points therefore controls the global dynamics, and characterizing it is a natural route to understanding the system. For deterministic dynamics, this local analysis is standard---the Jacobian's eigenvectors give normal modes and its eigenvalues give growth rates. For \emph{stochastic} dynamics one would like the analogous local theory of the slow, collective modes that govern how probability is transported through the neighborhood.

Those modes are the singular functions of the stochastic Koopman operator whose singular values lie close to one, the spectral objects underlying \emph{coherent sets}~\cite{froyland2013analytic,froyland2015dynamic,klus2019kernel}: subsets of phase space that remain dynamically coherent rather than dispersing. At present, coherent sets are studied almost exclusively \emph{numerically} and \emph{globally}: one discretizes a transfer or Koopman operator over the entire phase space of a nonlinear system---using Ulam-type discretizations~\cite{ulam1960collection,li1976finite}, finite elements~\cite{froyland2018robust}, extended dynamic mode decomposition~\cite{williams2015data}, or diffusion maps~\cite{coifman2006diffusion}---and reads off a partition~\cite{froyland2013analytic,klus2019kernel}. This yields neither a closed form nor a \emph{local law}: no statement of what the coherent modes look like near a given fixed point. The present paper supplies that missing local, analytic description.

A stable fixed point (a sink) linearizes to a contracting Ornstein--Uhlenbeck (OU) process: an isolated metastable well with no transport through it. The dynamically interesting fixed points are \emph{saddles}---hyperbolic points with at least one unstable direction. A saddle is a transition state: the bottleneck through which trajectories pass between basins, where reaction coordinates and committor functions are defined~\cite{hanggi1990reaction, vanden2006towards, bolhuis2002transition}. Its linearization is a continuous OU process whose drift matrix has eigenvalues of both signs, and this is the regime we solve.

The breakdown is easy to see. When the drift matrix has unstable eigenvalues, there is no normalizable stationary measure: probability escapes along the unstable manifold. As a consequence, the naive extension of singular functions from stable fixed points to saddles leads to non-integrable functions and singular values greater than one~\cite{pughe2026neurons} (Fig.~\ref{gamma}). This breaks the singular-value problem at its foundation. The singular value represents the correlation coefficient between the left singular function evaluated at the initial time and the right singular function evaluated at time $\tau$. In the coherent-set picture, where these functions approximate indicators of regions, this correlation is essentially the fraction of trajectories that begin in the coherent set and remain there after time $\tau$---a survival probability, which cannot exceed one. The breakdown is not merely analytic: the data-driven route fails in tandem, because the empirical covariances that any CCA-type estimator~\cite{hotelling1992relations,wu2020variational,klus2019kernel} relies on diverge along the unstable directions as the trajectories escape, so the correlation matrix has no finite limit and the canonical problem is ill-posed on the raw data (Fig.~\ref{gamma}). Any correct theory of the saddle case must first restore a proper Hilbert space.

Our resolution, and the conceptual core of the paper, is to \emph{condition} the dynamics: restricting to trajectories that have not yet escaped is a Doob's $h$-transform~\cite{doob1957conditional} that stabilizes the drift and makes the slow modes well defined. We show (i) that the required conditioning is the unique \emph{stabilizing} solution of an algebraic Riccati equation~\cite{lancaster1995algebraic}, selected automatically as the attractor of a Riccati flow (Sec.~\ref{sec:stab}); (ii) that the conditioned slow modes admit closed form whose decay rate is the escape rate along the unstable manifold, and that multiplication by $h$ is a unitary map relating them to ordinary Hermite modes (Sec.~\ref{sec:modes}); (iii) that this unitary structure yields a data-driven estimator---canonical correlation analysis on conditioned trajectories---which recovers the analytic modes exactly (Sec.~\ref{sec:datadriven}); (iv) that in the special case with a single slow unstable mode, the leading predictive direction converges to the unstable left eigenvector, justifying the reaction coordinate in classical theories (Sec.~\ref{sec:singlemode}); and (v) that the construction establishes bridges to transition-path theory, stochastic thermodynamics, linear control, and diffusion-model time reversal (Sec.~\ref{sec:discussion}).


\section{Setup and main result}
\label{sec:setup}

\subsection{Local linearization of stochastic nonlinear dynamics}

Consider a stochastic nonlinear system on $\R^d$,
\begin{equation}
\frac{d\mathbf x(t)}{dt} = \mathbf b(\mathbf x) + \bm\xi(t),
\qquad \langle\bm\xi(t)\bm\xi(t')^{\!\top}\rangle = D\,\delta(t-t'),
\label{eq:nonlinear}
\end{equation}
where $\mathbf b:\R^d\to\R^d$ is the drift and $D\succ0$ is the constant diffusion matrix.
Let $\mathbf x_*$ be a fixed point, $\mathbf b(\mathbf x_*)=0$. Writing $A=\nabla\mathbf b(\mathbf x_*)$ for the Jacobian and shifting the origin to $\mathbf x_*$, the local dynamics is the OU process
\begin{equation}
\frac{d\mathbf x(t)}{dt} = A\,\mathbf x(t) + \bm\xi(t).
\label{eq:ou}
\end{equation}
We assume for clarity that the eigenvalues and eigenvectors of $A$ are real (the complex case is discussed in Appendix S4 D). Denote by $\lambda_i$ the eigenvalues of $A$, by $r_i$ the right eigenvectors ($A r_i=\lambda_i r_i$), and by $\ell_i$ the left eigenvectors ($A^{\!\top}\ell_i=\lambda_i\ell_i$), normalized biorthogonally, $\ell_i^{\!\top}r_j=\delta_{ij}$. We split the spectrum into a \emph{stable} part (eigenvalues $\lambda_i^-<0$, eigenvectors $r_i^-,\ell_i^-$) and an \emph{unstable} part (eigenvalues $\lambda_i^+>0$, eigenvectors $r_i^+,\ell_i^+$), collect the unstable left eigenvectors into the columns of $\Lp$ and the unstable eigenvalues into $\Lcal=\mathrm{diag}(\lambda_i^+)$, so that $A^{\!\top}\Lp=\Lp\Lcal$. A saddle is the case where both parts are nonempty.

\subsection{Fokker--Planck operator and the Koopman generator}

The density $p(\mathbf x,t)$ obeys the Fokker--Planck equation $\partial_t p=\Lf p$ with the forward generator
\begin{equation}
\Lf f = -\nabla\!\cdot(A\mathbf x\,f) + \tfrac12 D^{\alpha\beta}\partial_\alpha\partial_\beta f .
\label{eq:FP}
\end{equation}
Its adjoint with respect to the Euclidean inner product $\langle f,g\rangle=\int f g\,d\mathbf x$ is the backward Kolmogorov operator, which is exactly the generator of the finite-time \emph{stochastic Koopman operator}:
\begin{equation}
\Lb g = (A\mathbf x)\cdot\nabla g + \tfrac12 D^{\alpha\beta}\partial_\alpha\partial_\beta g,
\qquad \Kop_\tau = e^{\Lb\tau},
\label{eq:generator}
\end{equation}
where the Koopman operator $\Kop_\tau$ propagates observables $g:\R^d\to\R$ by the conditional expectation~\cite{klus2024dynamical}
\begin{equation}
(\Kop_\tau g)(\mathbf x) \coloneq \E\big[\,g(\mathbf x(\tau))\;\big|\;\mathbf x(0)=\mathbf x\,\big].
\label{eq:koopman}
\end{equation}
A stationary measure satisfying $\Lf\rho_0=0$ is given by $\rho_0(\mathbf x)\propto\exp(-\tfrac12\mathbf x^{\!\top}\Sigma^{-1}\mathbf x)$, with $\Sigma$ solving the Lyapunov equation
\begin{equation}
A\Sigma + \Sigma A^{\!\top} = -D .
\label{eq:lyap}
\end{equation}
Eq.~\eqref{eq:lyap} has a unique solution iff no two eigenvalues of $A$ sum to zero~\cite{horn2012matrix}, and $\Sigma\succ0$ is the covariance iff $A$ is Hurwitz (all $\lambda_i<0$). If $A$ has unstable eigenvalues, $\Sigma$ is indefinite (negative definite if all modes are unstable); then $\rho_0$ grows in the unstable directions and is \emph{not normalizable}. Although such a non-normalizable $\rho_0$ does not define a probability distribution on $\mathbb{R}^d$, it remains a well-defined local stationary profile for the linearized dynamics near the fixed point. This is the root of the obstruction described in Sec.~I.

\subsection{Singular functions of the Koopman operator}

Coherent modes are defined variationally. Given a reference measure $\mu$ with $\langle f,g\rangle_\mu\coloneq\int fg\,d\mu$ and $\|f\|_\mu\coloneq\sqrt{\langle f,f\rangle_\mu}$, the leading singular value and functions solve
\begin{equation}
\sigma_1 = \max_{f,g\in L^2(\mu)}\frac{\langle f,\Kop_\tau g\rangle_\mu}{\|f\|_\mu\,\|g\|_\mu},
\label{eq:variational}
\end{equation}
where the maximizers are the leading left singular function $f_1$ and right singular function $g_1$. Further modes maximize the same ratio subject to $\mu$-orthogonality to the previous ones. Introducing the $\mu$-adjoint $\Kop_\tau^{\dagger}$ via $\langle f,\Kop_\tau g\rangle_\mu=\langle\Kop_\tau^{\dagger}f,g\rangle_\mu$, and the forward--backward and backward--forward operators $\mathcal F_\tau=\Kop_\tau\Kop_\tau^{\dagger}$ and $\mathcal B_\tau=\Kop_\tau^{\dagger}\Kop_\tau$, the singular functions are equivalently given by the following eigenfunctions~\cite{klus2024dynamical}:
\begin{equation}
\mathcal F_\tau f_i = \sigma_i^2 f_i,\qquad
\mathcal B_\tau g_i = \sigma_i^2 g_i,\qquad
\Kop_\tau g_i = \sigma_i f_i.
\label{eq:svd}
\end{equation}
Typically, $\mu$ is taken as the invariant measure. However, when $A$ has unstable modes, the invariant measure is non-normalizable, the naive stable-OU polynomial singular functions are not integrable with respect to this measure, and the ratio in Eq.~\eqref{eq:variational} can exceed one---the $\sigma>1$ pathology of Sec. I.

\subsection{Summary of the main result}
\label{sec:summary}

We show that conditioning replaces the ill-defined problem by a well-posed one on a stabilized process (Sec.~\ref{sec:stab}), whose singular functions are then obtained exactly. The result, derived in Sec.~\ref{sec:modes} and Appendix S8, is the following. Conditioning produces a stabilized drift $\Atil=A+DP_+$ (Hurwitz) with $P_+=\Lp(\Lp^{\!\top}\Sigma\Lp)^{-1}\Lp^{\!\top}$, a stabilized covariance $\Stil$ satisfying $\Atil\Stil+\Stil\Atil^{\!\top}=-D$, a Gaussian factor $h(\mathbf x)=\exp(\tfrac12\mathbf x^{\!\top}P_+\mathbf x)$, and canonical directions $v_i,u_i$ with canonical correlations $\gamma_i\in(0,1)$ (the singular values of the stabilized process, defined in Sec.~\ref{sec:cca}). We introduce a reference weight $\rhotil_0$ defined in Sec.~\ref{sec:rhotil}, the natural extension of $\rho_0$ to saddles, on which $\Kop_\tau$ is a contraction. The saddle Koopman singular functions and singular values, with respect to $\rhotil_0$, are indexed by a multi-index $\bm n\in\mathbb N_0^{\,d}$ and given in closed form by
\begin{equation}
\boxed{\;
\begin{aligned}
\Phi^{L}_{\bm n}(\mathbf x) &= \frac{h(\mathbf x)}{\sqrt{\bm n!}}\prod_{i}\He_{n_i}\!\big(v_i^{\!\top}\mathbf x\big),\\[2pt]
\Phi^{R}_{\bm n}(\mathbf x) &= \frac{h(\mathbf x)}{\sqrt{\bm n!}}\prod_{i}\He_{n_i}\!\big(u_i^{\!\top}\mathbf x\big),\\[2pt]
\sigma_{\bm n} &= e^{-\Tr(\Lcal)\,\tau}\prod_{i}\gamma_i^{\,n_i},\\[2pt]
e^{\Atil^{\!\top}\tau}&\Stil^{-1}e^{\Atil\tau}\Stil v_i =\gamma_i^2 v_i,\\[2pt]
\Stil^{-1}&e^{\Atil\tau}\Stil\,e^{\Atil^{\!\top}\tau}u_i =\gamma_i^2 u_i,
\end{aligned}
\;}
\label{eq:mainresult}
\end{equation}
where $\He_n$ is the probabilists' Hermite polynomial and $\bm n!=\prod_i n_i!$. The Gaussian factor $h$ decays along the unstable manifold, with the decaying prefactor of the singular value $e^{-\Tr(\Lcal)\tau}$ set by the escape rate $\Tr(\Lcal)=\sum_i\lambda_i^+$. The remainder of the paper constructs each ingredient and proves Eq.~\eqref{eq:mainresult}.

\section{Stabilization by conditioning}
\label{sec:stab}

\subsection{Doob's \texorpdfstring{$h$}{h}-transform}
\label{sec:doob}

\emph{Conditioning on a single endpoint.} Retain only trajectories arriving at a prescribed point $\mathbf x(T)=\mathbf y$. The conditioned process is again a diffusion---Doob's $h$-transform---with the same $D$ and a drift modified by the gradient of the log transition density~\cite{doob1957conditional}. Writing $h(\mathbf x,t)=p(\mathbf y,T\mid\mathbf x,t)$ for the probability of reaching $\mathbf y$ from $(\mathbf x,t)$, the drift is modified as
\begin{equation}
A\mathbf x \;\longmapsto\; A\mathbf x + D\,\nabla\log h(\mathbf x,t).
\label{eq:driftshift-point}
\end{equation}

\emph{Conditioning on a final distribution.} More generally, one can condition on a prescribed terminal distribution $\pi(\mathbf x_T)$. The $h$-function in this case is given by the average of the single-point $h$-function weighted by the conditional-to-original distribution ratio (Appendix S1),
\begin{equation}
h(\mathbf x,t) = \int p(\mathbf x_T,T\mid\mathbf x,t)\,\frac{\pi(\mathbf x_T)}{p(\mathbf x_T)}\,d\mathbf x_T
= (\Kop_{T-t}\,w)(\mathbf x),
\label{eq:hfinal}
\end{equation}
where $w\coloneq \frac{\pi}{p}$ is the reweighting coefficient. Thus, the conditioning $h$-function is itself a Koopman image of the reweighting $w$. For a Gaussian $w(\mathbf x)\propto\exp(-\tfrac12\mathbf x^{\!\top}\Sigma'^{-1}\mathbf x)$, the long-horizon limit, $T\gg t$, yields a Gaussian $h$,
\begin{equation}
\begin{aligned}
&h(\mathbf x, t) \propto h(\mathbf x)=\exp\!\big(\tfrac12\mathbf x^{\!\top}P_+\mathbf x\big),\\
&P_+ = \Lp\big(\Lp^{\!\top}\Sigma\Lp\big)^{-1}\Lp^{\!\top},
\end{aligned}
\label{eq:hP}
\end{equation}
with $\Lp$ the unstable left eigenvectors of Sec.~\ref{sec:setup}. The proof and explanation are given in the next subsection. We stress that this long-horizon limit is independent of the reweighting window $\Sigma'$. Following Eq. \eqref{eq:driftshift-point}, the conditioned drift matrix is modified as
\begin{equation}
A \;\longmapsto\; \Atil = A + DP_+.
\label{eq:Atil}
\end{equation}
The two extremes are transparent:
\begin{align}
\text{all stable:}\quad & \Lp=\varnothing,\quad P_+=0,\quad \Atil=A, \label{eq:corner-stable}\\
\text{all unstable:}\quad & P_+=\Sigma^{-1},\quad \Atil=-\Sigma A^{\!\top}\Sigma^{-1}. 
\label{eq:corner-unstable}
\end{align}
In the all-stable case, conditioning does nothing. In the all-unstable case, the conditioned drift is exactly the score-based reverse-time SDE of generative diffusion modeling~\cite{anderson1982reverse,song2020score}: $\Atil \mathbf x=A\mathbf x-D\nabla\log\rho_0(\mathbf{x})$. The saddle interpolates between them.

\begin{figure}[t]
  \centering
  \begin{subfigure}[t]{0.53\linewidth}
    \caption{}
    \includegraphics[width=\linewidth]{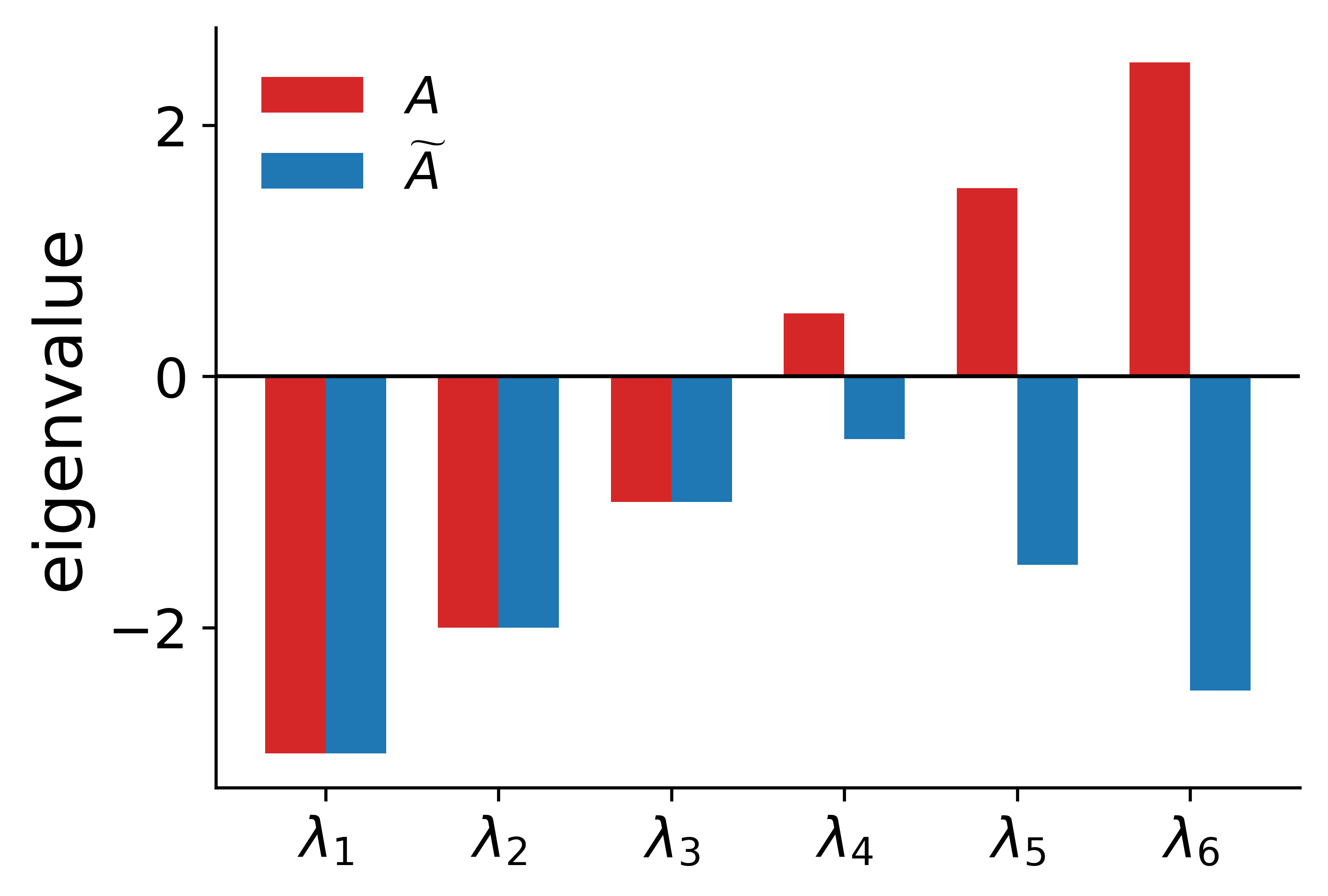}
    \label{eigenvalue reversal}
    \end{subfigure}
    \begin{subfigure}[t]{0.45\linewidth}
    \caption{}
    \includegraphics[width=\linewidth]{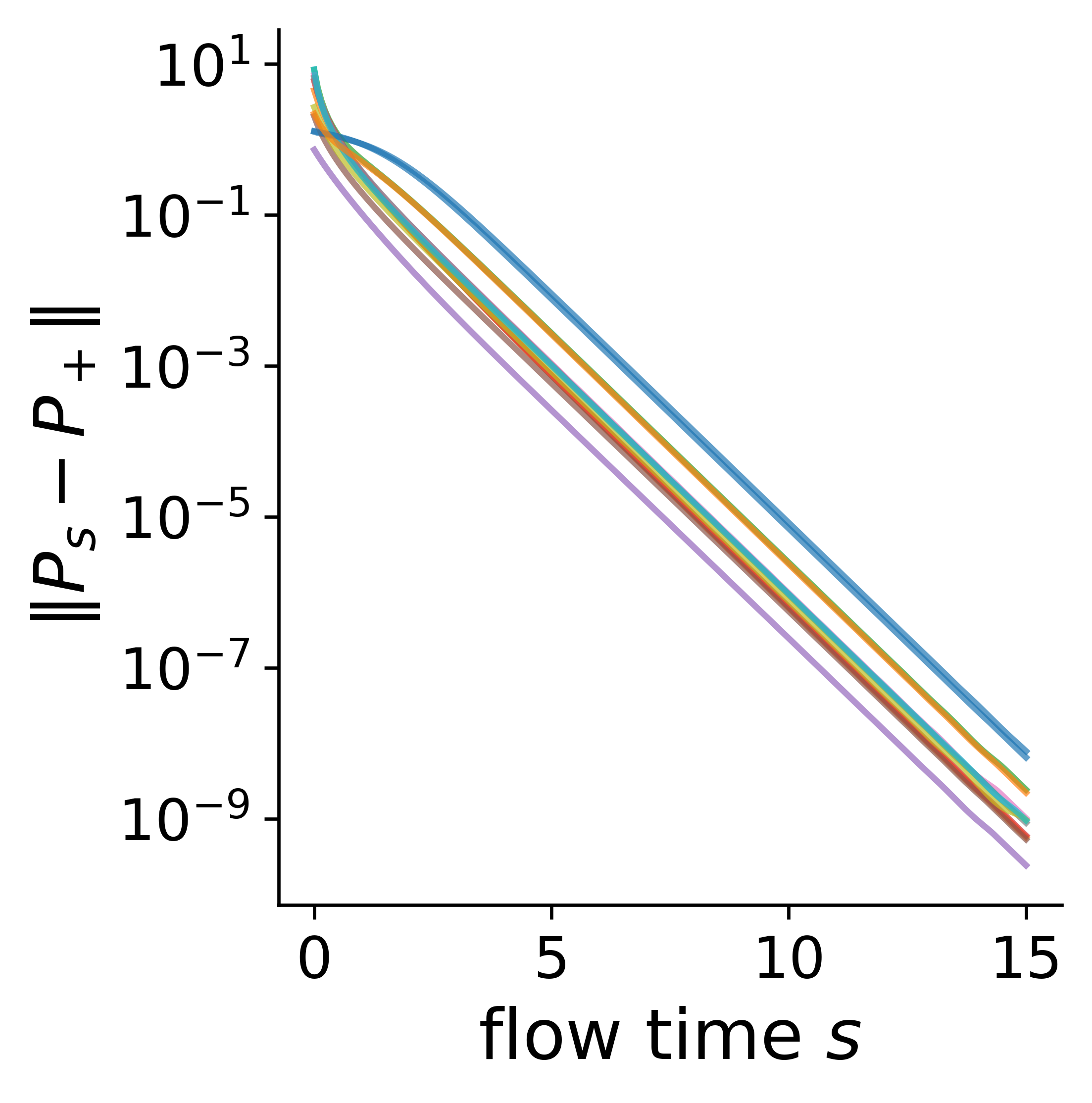}
    \label{reccati flow}
    \end{subfigure}
  \caption{Stabilization through Riccati flow. (a) Negative eigenvalues of $A$ are kept in $\Atil$, while positive eigenvalues are flipped. (b) The Frobenius norm of the difference between $P_s$ and the attractor $P_+$ (Eq.~\eqref{eq:hP}) decays in the Riccati flow. Different traces correspond to different initial conditions $\Sigma'$.}
  \label{stabilization}
\end{figure}

\begin{proposition}[Selective time reversal]
\label{prop:reflect}
With $\Atil$ from Eq.~\eqref{eq:Atil},
\begin{equation}
\Atil\,r_i^- = \lambda_i^-\,r_i^-,
\qquad
\Atil\,(\Sigma\ell_i^+) = -\lambda_i^+\,(\Sigma\ell_i^+):
\end{equation}
the stable eigenvalues are unchanged and the unstable ones are reflected across zero. In particular, $\Atil$ is Hurwitz.
\end{proposition}
\noindent Proof in Appendix S3. Conditioning is thus a \emph{selective time reversal}, flipping the arrow of time only along the escaping directions. The numerical verification is given in Fig.~\ref{eigenvalue reversal}, and Fig.~\ref{linear demo} shows the original and stabilized drift of an example linear saddle system.

\subsection{The Riccati flow and its stabilizing attractor}
\label{sec:riccati}

We now show why the conditioning selects precisely $P_+$ in Eq.~\eqref{eq:hP}. From Eq.~\eqref{eq:hfinal}, $h_s= \Kop_s w$ with $s=T-t$; the OU semigroup preserves the Gaussian class, so $h_s\propto\exp(\tfrac12\mathbf x^{\!\top}P_s\mathbf x)$ with $P_s$ symmetric, and we have (Appendix S2)
\begin{equation}
    \Lb h_s = \Big[\tfrac12 x^{\!\top}\!\big(A^{\!\top}P_s+P_sA+P_sDP_s\big)x+\tfrac12\Tr(DP_s)\Big]h_s.
\end{equation}
Substituting into $\partial_s h_s=\Lb h_s$ gives a matrix \emph{Riccati flow} for the quadratic form (Appendix S5),
\begin{equation}
\dot P_s = A^{\!\top}P_s + P_s A + P_s D P_s .
\label{eq:riccatiflow}
\end{equation}
As $s\to\infty$ this flow converges to a fixed point $P$ solving the \emph{algebraic Riccati equation}~\cite{lancaster1995algebraic}
\begin{equation}
A^{\!\top}P + PA + PDP = 0,
\label{eq:riccati}
\end{equation}
and thus $\exp\!\big(\tfrac12\mathbf x^{\!\top}P\mathbf x\big)$ is an eigenfunction of $\Lb$ with eigenvalue $\frac{1}{2}\Tr(DP)$.

To see which fixed point is reached, we first describe the structure of \emph{all} solutions, then identify the attractor.

\emph{Structure of the solutions.} Since $P$ enters through the quadratic form, we restrict throughout to symmetric solutions $P=P^\top$. Each symmetric solution $P$ of Eq.~\eqref{eq:riccati} produces a conditioned drift $\Atil=A+DP$ whose eigenvalues have a simple relation to those of $A$: every eigenvalue $\lambda_i$ of $A$ appears in $\Atil$ either unchanged (``kept'') or negated to $-\lambda_i$ (``flipped''). The one-dimensional case makes this concrete: for $\dot x=\lambda x+\xi$ the Riccati equation $2\lambda P+DP^2=0$ has the two roots $P=0$ (giving $\Atil=\lambda$, kept) and $P=-2\lambda/D=\Sigma^{-1}$ (giving $\Atil=-\lambda$, flipped). In $d$ dimensions, a solution is specified by choosing, independently for each of the $d$ modes, whether to keep or flip it. There are therefore $2^d$ solutions, given explicitly by
\begin{equation}
P_{\mathcal F} = L_{\mathcal F}\big(L_{\mathcal F}^{\!\top}\Sigma L_{\mathcal F}\big)^{-1}L_{\mathcal F}^{\!\top},
\label{eq:Pfamily}
\end{equation}
where $\mathcal F\subseteq\{1,\dots,d\}$ is the set of flipped modes and $L_{\mathcal F}$ collects their left eigenvectors. Three landmarks: $\mathcal F=\varnothing$ gives $P=0$; the full set gives $P=\Sigma^{-1}$ (so $\Atil=-\Sigma A^{\!\top}\Sigma^{-1}$); and $\mathcal F=\{\text{unstable modes}\}$ gives $P_+$ in Eq.~\eqref{eq:hP}. The systematic derivation via the Hamiltonian matrix $\bigl[\begin{smallmatrix}A&D\\0&-A^{\!\top}\end{smallmatrix}\bigr]$ and its invariant subspaces is in Appendix S4~\cite{lancaster1995algebraic}.

\emph{Which one does the flow reach?} A conditioned process that does not escape must be \emph{contracting}, i.e.\ $\Atil$ must be Hurwitz (all eigenvalues negative). Among the $2^d$ choices, exactly one achieves this: flip every unstable mode ($\lambda_i^+>0\mapsto-\lambda_i^+<0$) and keep every stable mode. This is the \emph{stabilizing} solution, Eq.~\eqref{eq:hP}. The mathematical reason for the Riccati flow~\eqref{eq:riccatiflow} converging to it and to no other is the following. Linearizing the flow about any solution $P_{\mathcal F}$ gives growth rates $\widetilde\lambda_i+\widetilde\lambda_j$, where $\{\widetilde\lambda_i\}$ are the eigenvalues of $\Atil$, and these are all negative---so $P_{\mathcal F}$ is attracting---only when every $\widetilde\lambda_i<0$, i.e.\ only for the stabilizing choice; every other solution has a growing direction and is a saddle of the flow (Appendix S5). Hence, long-horizon conditioning lands on Eq.~\eqref{eq:hP} regardless of the conditioning window $\Sigma'$. We numerically verify the convergence for different windows $\Sigma'$, as shown in Fig.~\ref{reccati flow}.

\subsection{The reference measure \texorpdfstring{$\rhotil_0$}{rhotil0}}
\label{sec:rhotil}

The stabilized process $\dot{\mathbf x}=\Atil\mathbf x+\bm\xi$ has a stationary Gaussian density $\rho_{\Atil}\propto\exp(-\tfrac12\mathbf x^{\!\top}\Stil^{-1}\mathbf x)$, $\Atil\Stil+\Stil\Atil^{\!\top}=-D$, $\Stil\succ0$. In the all-stable case, $\Stil=\Sigma$; in the all-unstable case, $\Stil=-\Sigma$. As illustrated by the following proposition, $\Stil$ can generically be expressed in closed form, through the pair of complementary Riccati solutions.
\begin{proposition}[Closed form for the stabilized covariance]
\label{prop:Stil}
\begin{equation}
\Stil=(P_--P_+)^{-1},
\label{eq:Stilinv}
\end{equation}
where $P_\pm=L_\pm\big(L_\pm^{\!\top}\Sigma L_\pm\big)^{-1}L_\pm^{\!\top}$, with $L_\pm$ the collection of unstable/stable left eigenvectors of $A$.
\end{proposition}
\noindent Proof in Appendix S6. In contrast to $P_+$, $P_-$ solely flips all the stable modes (Eq.~\eqref{eq:Pfamily} with $\mathcal F=\{\text{stable}\}$); correspondingly $A+DP_-$ is purely unstable. Eq.~\eqref{eq:Stilinv} is the covariance-level counterpart of the selective time reversal of Proposition~\ref{prop:reflect}. Thus, $\rho_{\Atil}\propto \exp\!\Big[-\tfrac12\mathbf x^{\!\top}\big(P_--P_+\big)\mathbf x\Big]$.

We now consider the reference measure defined as
\begin{equation}
\begin{aligned}
\rhotil_0(\mathbf x) =& \rho_{\Atil}(\mathbf x)\,h(\mathbf x)^{-2}\\
=& \exp\!\big(-\tfrac12\mathbf x^{\!\top}\Stil^{-1}\mathbf x\big)\,\exp\!\big(-\mathbf x^{\!\top}P_+\mathbf x\big)\\
=& \exp\!\Big[-\tfrac12\mathbf x^{\!\top}\big(P_-+P_+\big)\mathbf x\Big].\\
\end{aligned}
\label{eq:rhotil}
\end{equation}
The all-stable case gives $P_-=\Sigma^{-1},\: P_+=0$ and the all-unstable case gives $P_-=0,\: P_+=\Sigma^{-1}$. Thus, in both all-stable and all-unstable extremes, $\rhotil_0=\rho_0$. In the genuine saddle case, they are generically not identical because the stable and unstable directions demand opposite curvatures, and $\rhotil_0=\rho_{\Atil}h^{-2}$ is the minimal deformation that honestly keeps the stable manifold's Gaussian and the unstable manifold's reciprocal Gaussian. Thus, $\rhotil_0$ is the natural extension of $\rho_0$ from stable/unstable extremes to saddles, which decouples the stable and unstable covariances. We restrict the variational maximization for singular functions, Eq.~\eqref{eq:variational}, to the Hilbert space $L^2(\rhotil_0)$. The factorization Eq.~\eqref{eq:rhotil} makes the relevant integrability transparent: for any polynomial $\mathcal H$,
\begin{equation}
\|h\,\mathcal H\|^2_{L^2(\rhotil_0)}=\int h^2|\mathcal H|^2\,\rhotil_0\,d\mathbf x
=\int|\mathcal H|^2\,\rho_{\Atil}\,d\mathbf x<\infty.
\end{equation}
In particular, the singular functions in Eq.~\eqref{eq:mainresult} are square-integrable in $L^2(\rhotil_0)$. As shown next, this specific measure is exactly the one making conditioning a unitary change of representation, and it is what the data-driven method computes.

\section{Exact singular modes and their unitary origin}
\label{sec:modes}

The main result is Eq.~\eqref{eq:mainresult}: the saddle Koopman singular functions with respect to $\rhotil_0$ are products of univariate Hermite polynomials in the canonical coordinates $v_i^{\!\top}\mathbf x$, $u_i^{\!\top}\mathbf x$, multiplied by the Gaussian factor $h$, with singular values $e^{-\Tr(\Lcal)\tau}\prod_i\gamma_i^{n_i}$. A direct derivation is given in Appendix S8. Here we explain the result through the structure that makes it true, which also motivates the estimator of Sec.~\ref{sec:datadriven}.

As shown in Sec.~\ref{sec:riccati}, $h$ is a positive eigenfunction of the Koopman generator, for $P_+$ solving Eq.~\eqref{eq:riccati},
\begin{equation}
\Lb h = \lambda\,h,\qquad \lambda = \tfrac12\Tr(DP_+) = -\Tr(\Lcal).
\label{eq:heig}
\end{equation}
It follows that the Koopman generator of the stabilized process is given by
\begin{equation}
\begin{aligned}
    \widetilde{\Lf}^\dagger g=&\frac{1}{h}\Lb(hg)-\frac{1}{h}\Lb(h)g\\
    =& \frac{1}{h}\Lb(hg)+\Tr(\Lambda_+)g,
\end{aligned}
\end{equation}
which leads to the following proposition.
\begin{proposition}[Conditioning is a unitary change of representation]
\label{prop:unitary}
Let $M_h f=hf$. Because $\rhotil_0=\rho_{\Atil}h^{-2}$,
$M_h:\ L^2(\rho_{\Atil})\to L^2(\rhotil_0)\ $ is unitary:
\begin{equation}
\|M_h f\|_{L^2(\rhotil_0)}=\|f\|_{L^2(\rho_{\Atil})},
\end{equation}
and it intertwines the Koopman operators of the original and stabilized processes,
\begin{equation}
\Kop_\tau = e^{-\Tr(\Lcal)\tau}\,M_h\,\widetilde{\Kop}_\tau\,M_h^{-1}.
\label{eq:intertwine}
\end{equation}
\end{proposition}
\noindent Proof in Appendix S7. Since the singular value decomposition is invariant under unitary conjugation and a scalar prefactor merely rescales singular values, the singular functions of the saddle Koopman operator (with respect to $\rhotil_0$) are $M_h$ applied to those of the \emph{stabilized} process (with respect to $\rho_{\Atil}$), with singular values scaled by $e^{-\Tr(\Lcal)\tau}$. It also follows from the unitarity that $\Kop_\tau$ is a contraction on $L^2(\rhotil_0)$, so the singular value problem is well-posed (Appendix S9).

The stabilized problem is a stable-OU singular-value problem, and its singular functions are products of univariate Hermite polynomials~\cite{pughe2026neurons}. The first-order left and right singular functions are the canonical scores $v_i^{\!\top}\mathbf x$ and $u_i^{\!\top}\mathbf x$, respectively, where $v_i,u_i$ solve
\begin{equation}
\begin{aligned}
&F v_i = \gamma_i^2 v_i,\quad B u_i = \gamma_i^2 u_i,\\
&F =e^{\Atil^{\!\top}\tau}\Stil^{-1}e^{\Atil\tau}\Stil,\quad
B = \Stil^{-1}e^{\Atil\tau}\Stil\,e^{\Atil^{\!\top}\tau},
\end{aligned}
\label{eq:FB}
\end{equation}
and $\gamma_i\in(0,1)$ are the canonical correlations in descending order $\gamma_1\ge\gamma_2\ge\cdots$. Crucially, in the canonical coordinates, the stationary measure factorizes into independent standard Gaussians, so the higher-order singular functions are simply \emph{products} of univariate Hermite polynomials (proof in Appendix S8).

\emph{The leading coherent mode is linear.} Because $\gamma_i\in(0,1)$, excluding the trivial mode $\bm n=\bm 0$, the product $\prod_i\gamma_i^{n_i}$ is largest only if there is a single mode, i.e., $\sum_i n_i=1$. That is, the largest nontrivial singular value is $\sigma_{\bm e_1}=e^{-\Tr(\Lcal)\tau}\gamma_1$, attained by a \emph{first-order} mode: $\Phi^L_1=h\,v_1^{\!\top}\mathbf x$. Higher-order modes are strictly less coherent. Because coherent sets are partitioned by the sign pattern of top singular functions, and $h>0$, this explains why the linear filter---and hence canonical correlation analysis, which computes exactly the first-order modes (Sec.~\ref{sec:datadriven})---suffices to identify the dominant coherent partition.

\section{Data-driven estimator}
\label{sec:datadriven}

\subsection{Koopman singular functions and canonical correlation analysis}
\label{sec:cca}

The mainstream data-driven methods for approximating transfer operators from trajectory data---extended dynamic mode decomposition (EDMD)~\cite{williams2015data}, the variational approach for Markov processes (VAMP)~\cite{wu2020variational}, and other variants---are, at their core, correlation methods: they lift the state through a dictionary of observables and solve a generalized eigenproblem of empirical instantaneous and time-lagged covariance matrices. When the object sought are the leading singular functions rather than eigenfunctions, this reduction is exactly a canonical correlation analysis (CCA) between the observables at the initial and final times~\cite{wu2020variational,klus2019kernel}. Our estimator is of the same family. The point of this section is that, for the saddle case, the correct CCA is performed on the \emph{conditioned} (stabilized) trajectories rather than the raw ones. We first make the CCA--Koopman correspondence precise.

CCA~\cite{hotelling1992relations} finds linear projections of two random vectors that are maximally correlated. Applied in time to $\mathbf x(0)$ and $\mathbf x(\tau)$ of a stationary process, it seeks $v,u$ maximizing
\begin{equation}
\begin{aligned}
&\gamma_1 = \max_{v,u}\frac{v^{\!\top}C_\tau u}{\sqrt{(v^{\!\top}C_0 v)(u^{\!\top}C_0 u)}},\\
&C_0=\E[\mathbf x(0)\mathbf x(0)^{\!\top}],\ \ C_\tau=\E[\mathbf x(0)\mathbf x(\tau)^{\!\top}],
\end{aligned}
\end{equation}
with the canonical directions and the canonical correlations given by the eigenproblems of $C_0^{-1}C_\tau C_0^{-1}C_{-\tau}$ and $C_0^{-1}C_{-\tau} C_0^{-1}C_{\tau}$. For a stable OU process, this is exactly the first-order singular-value problem, Eqs.~\eqref{eq:svd} and~\eqref{eq:FB}, with $\mu=\rho_0$: substituting $C_0=\Sigma$, $C_\tau=\Sigma e^{A^{\!\top}\tau}$, $C_{-\tau}=e^{A\tau}\Sigma$,
\begin{equation}
\begin{aligned}
&C_0^{-1}C_\tau C_0^{-1}C_{-\tau} = e^{A^{\!\top}\tau}\Sigma^{-1}e^{A\tau}\Sigma = F,\\
&C_0^{-1}C_{-\tau} C_0^{-1}C_{\tau} = \Sigma^{-1}e^{A\tau}\Sigma e^{A^{\!\top}\tau} = B.
\end{aligned}
\label{eq:cca-equals-F}
\end{equation}
Thus, the CCA canonical scores $v^{\!\top}\mathbf x$, $u^{\!\top}\mathbf x$ are the first-order Koopman singular functions and the canonical correlations $\gamma_i$ are the singular values~\cite{akaike1975markovian,wu2020variational,pughe2026neurons}. In short, CCA is the empirical singular value decomposition of the Koopman operator for stable OU processes, which is why a data-driven procedure can recover the analytic modes. However, one cannot naively extend such CCA on raw trajectories to saddle cases, because the correlation matrices, $C_0$ and $C_\tau$, diverge along unstable directions (Fig.~\ref{gamma}).

\subsection{Algorithm for coherent-set partition}

In this subsection, we present the correct data-driven algorithm for estimating the coherent-set partition in the local neighborhood of a fixed point. Coherent sets are partitioned by the sign pattern of top Koopman singular functions~\cite{froyland2013analytic}. The estimation of canonical coordinates, $\xi_i=v_i^{\!\top}\mathbf x$ and $\eta_i=u_i^{\!\top}\mathbf x$ in Eq.~\eqref{eq:mainresult}, is all that the data-driven method needs: every singular function is a product of Hermite polynomials in them, $\Phi^{L}_{\bm n}\propto h\prod_i\He_{n_i}(\xi_i)$ and $\Phi^{R}_{\bm n}\propto h\prod_i\He_{n_i}(\eta_i)$. Because the Gaussian factor $h>0$ never changes sign, coherent sets are the sign sets of the Hermite part alone, i.e., they are determined entirely by the estimated canonical coordinates. Thus, once the $\{(v_i,u_i)\}$ are estimated, any coherent-set partition is available in closed form without further fitting. According to Eq.~\eqref{eq:FB}, the estimation of $\{(v_i,u_i)\}$ is a linear CCA on the stabilized process, obtained from the original (unstable) trajectories by the conditioning reweighting, $w$ in Eq.~\eqref{eq:hfinal}. This gives the following estimator.

\begin{center}
\fbox{\parbox{0.94\linewidth}{
\textbf{Algorithm: CCA with conditioning.}
\begin{enumerate}\itemsep2pt
\item Take the fixed point of interest as the center.
\item Collect trajectories $\{\mathbf x^{(m)}(t)\}_{t\in[0,T]}$ from any normalizable initial distribution.
\item Choose a confining covariance $\Sigma'$; assign trajectory $m$ the weight $w_m\propto\mathcal N\!\big(\mathbf x^{(m)}(T);0,\Sigma'\big)$, implementing conditioning on non-escape.
\item Choose lags $0\ll t_1<t_2\ll T$. Under the reweighted empirical measure $\{w_m\}$, form the weighted instantaneous and time-lagged covariances of $(\mathbf x(t_1),\mathbf x(t_2))$ and solve the CCA eigenproblem for the leading canonical pairs $\{(v_i,u_i)\}$.
\item Build any desired singular function from the estimated canonical coordinates via Eq.~\eqref{eq:mainresult}, and partition phase space by its sign pattern.
\end{enumerate}
}}
\end{center}

\begin{figure}[t]
  \centering
  \begin{subfigure}[t]{0.8\linewidth}
    \caption{}
    \includegraphics[width=\linewidth]{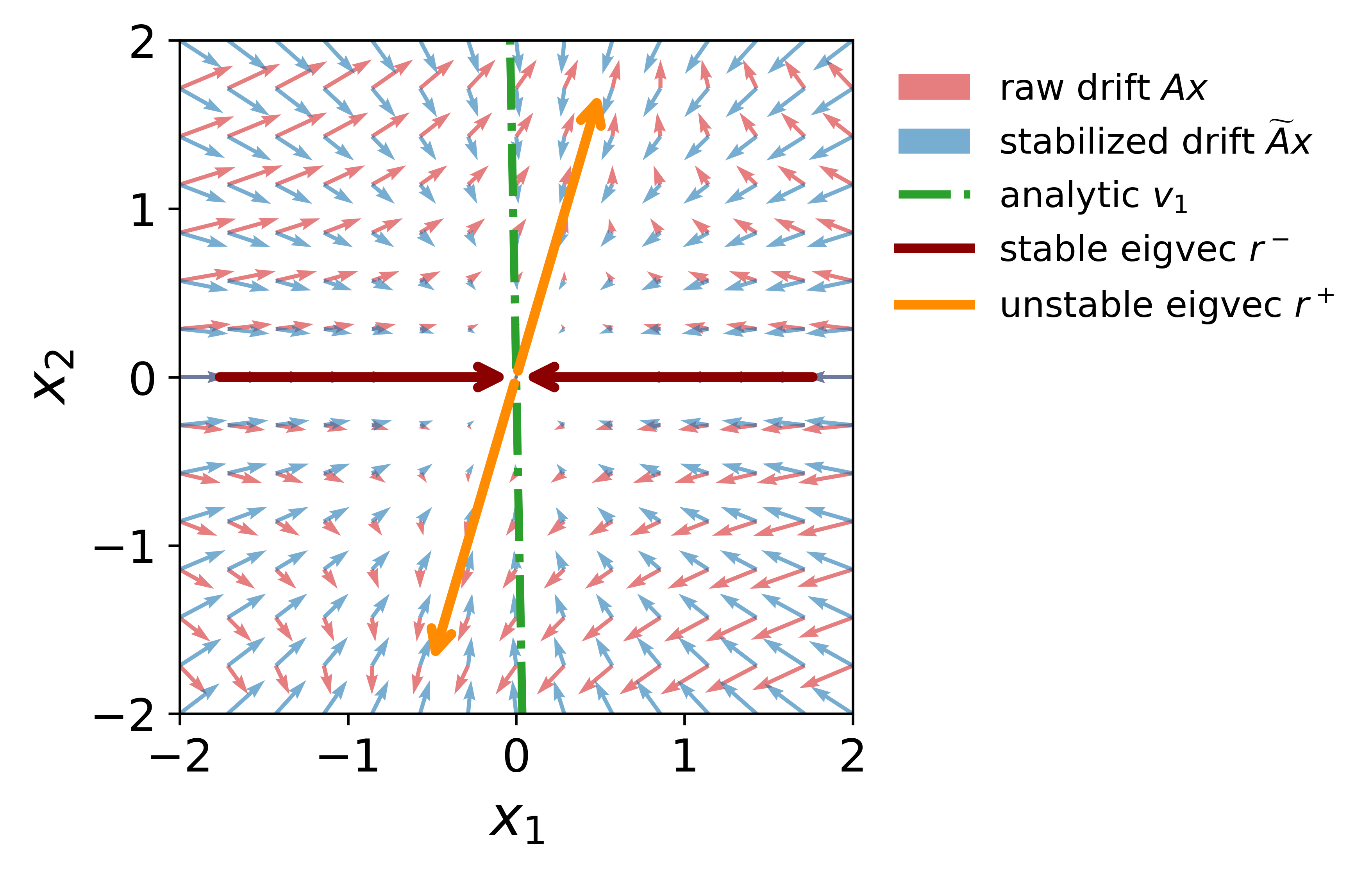}
    \label{linear demo}
    \end{subfigure}
    \begin{subfigure}[t]{0.51\linewidth}
    \caption{}
    \includegraphics[width=\linewidth]{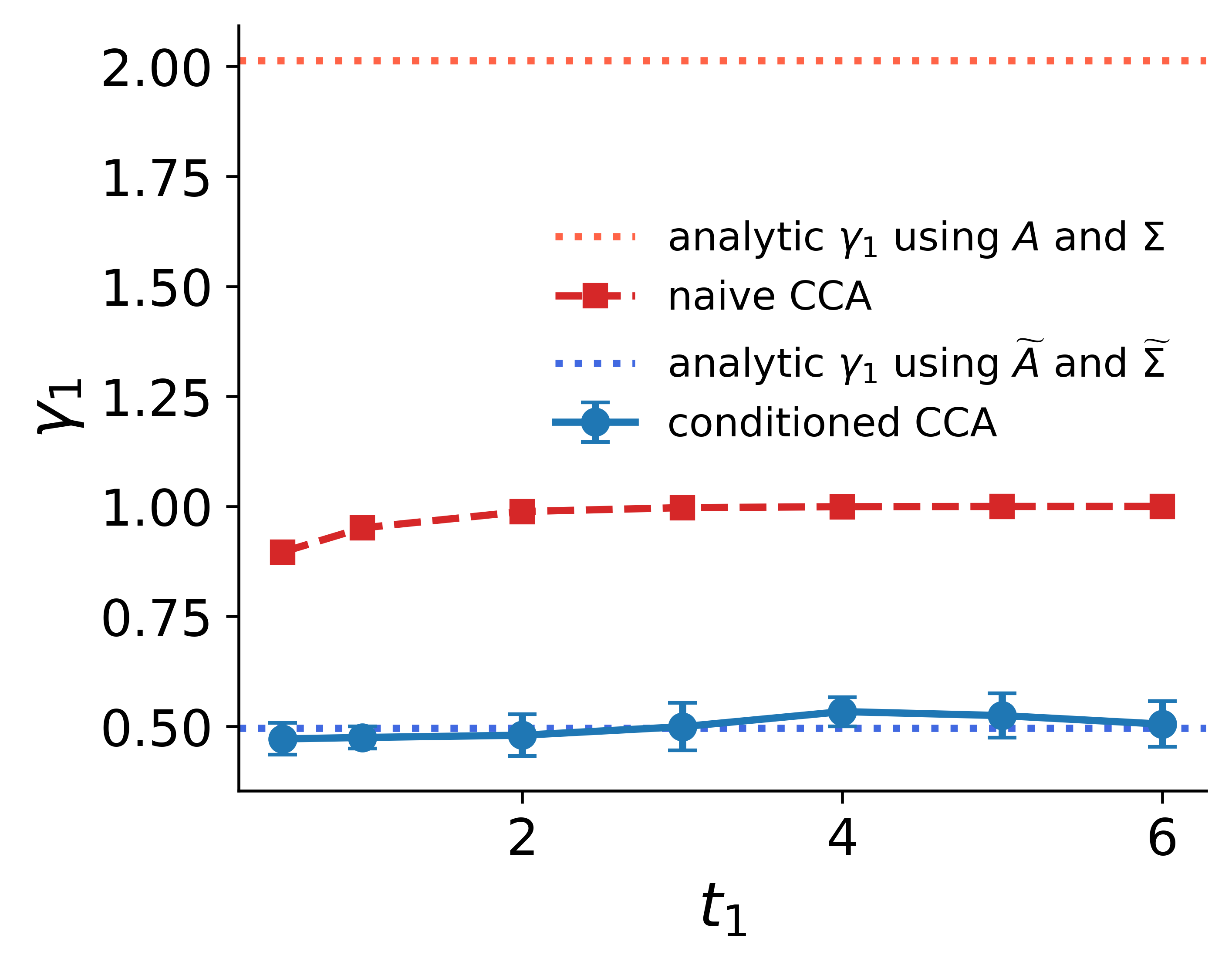}
    \label{gamma}
    \end{subfigure}
    \begin{subfigure}[t]{0.47\linewidth}
    \caption{}
    \includegraphics[width=\linewidth]{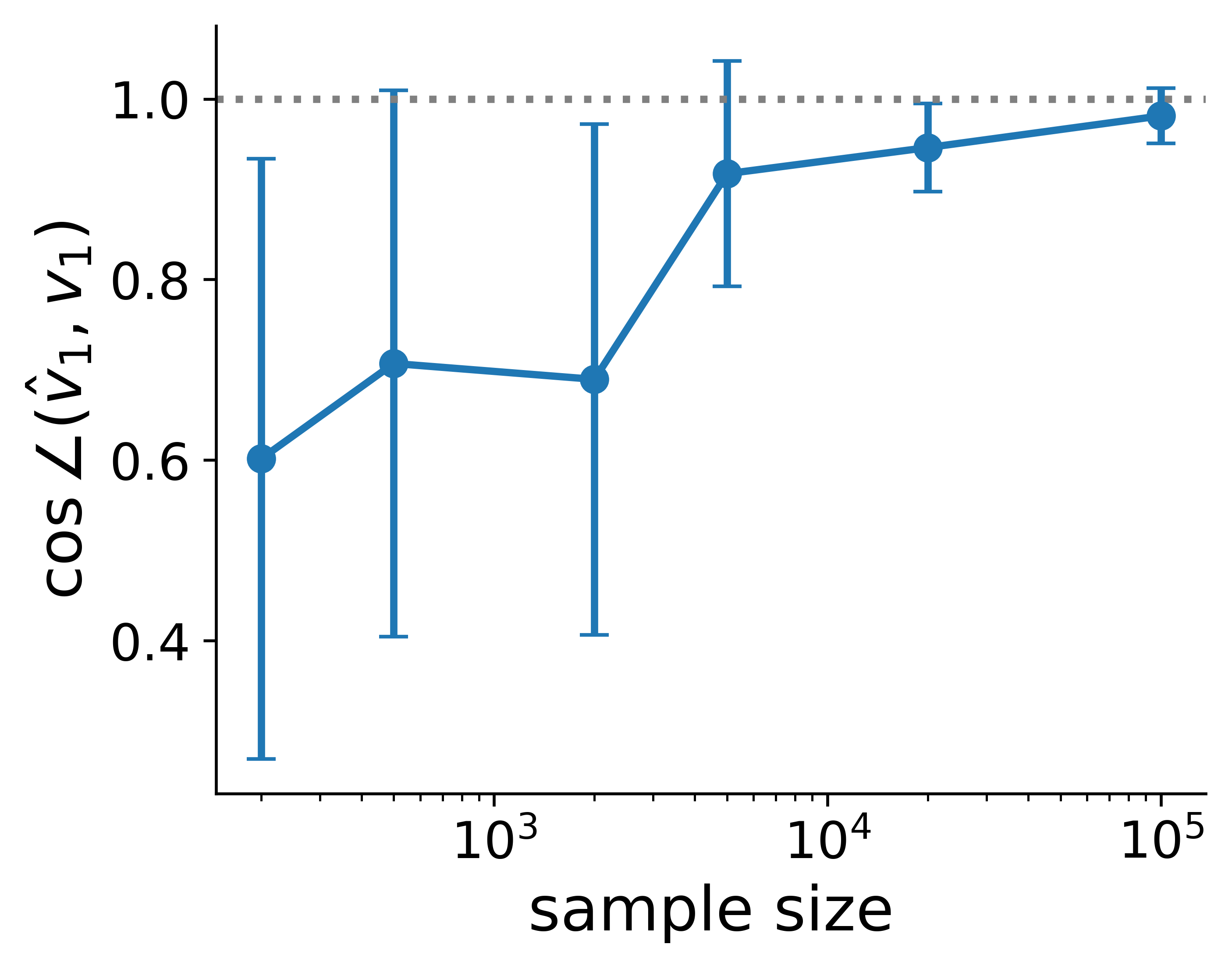}
    \label{linear cosine size}
    \end{subfigure}
  \caption{Linear stochastic dynamical system. (a) Demonstration of a 2D linear saddle system with $A=\begin{bmatrix}
-1. & 0.5 \\
0. & 0.7
\end{bmatrix}$ and 
$D=\begin{bmatrix}
1. & 0.3 \\
0.3 & 1.
\end{bmatrix}$. The eigenvectors, raw and stabilized drift, and analytic $v_1$ are plotted. (b) Canonical correlations calculated analytically (Eq.~\eqref{eq:FB}) using raw/stabilized dynamics and computed numerically using naive/conditioned CCA for varying $t_1$ and fixed $\tau$. (c) The cosine similarity between analytic $v_1$ and empirical $\hat{v}_1$ for increasing sample sizes.}
  \label{linear system}
\end{figure}

\begin{figure*}[t]
  \centering
  \includegraphics[width=0.98\linewidth]{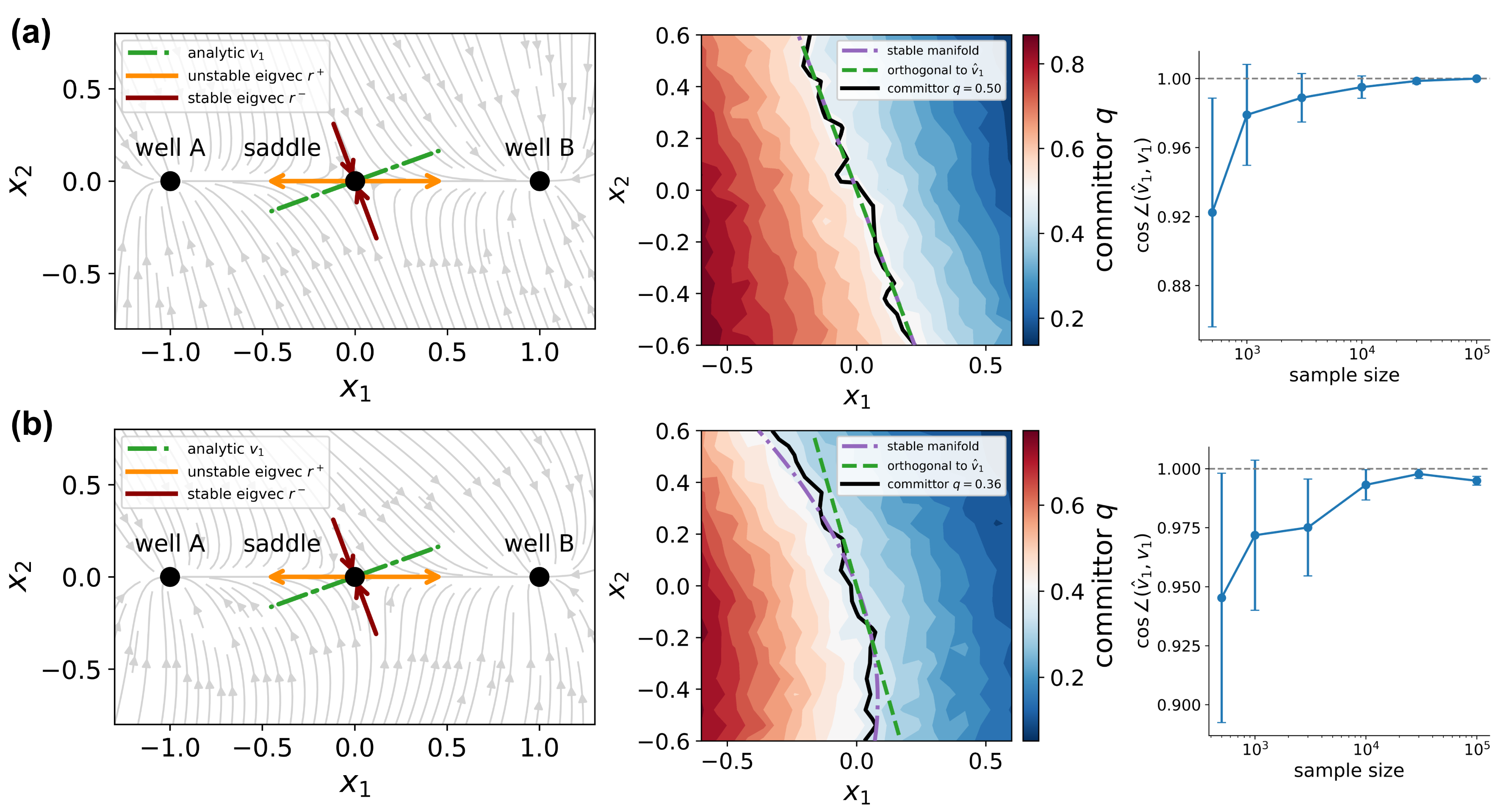}
  \caption{Nonlinear stochastic dynamical system. (a) $\frac{dx_1}{dt} = x_1 - x_1^3 - x_1x_2^2 + 1.5x_2$, $\frac{dx_2}{dt} = -3(1+x_1^2)x_2$. $D=I_2$. Left: a phase portrait of the dynamical system with two stable fixed points and a saddle. The eigenvectors and the analytic $v_1$ of the saddle are plotted. Middle: the heatmap shows the estimated committor $q(\mathbf{x})$ with the two stable fixed points as basins. The stable manifold, the committor contour passing the saddle, and the hyperplane orthogonal to $\hat{v}_1$ partitioning coherent sets are plotted. Right: The cosine similarity between analytic $v_1$ and empirical $\hat{v}_1$ of the saddle for increasing sample sizes. (b) Similar to (a) for $\frac{dx_1}{dt} = x_1 - x_1^3 + 1.5x_2 + 3x_2^2$, $\frac{dx_2}{dt} = -3(1+x_1^2)x_2$, where the stable manifold is curved.}
  \label{nonlinear system}
\end{figure*}

\emph{Leading linear partition.} The most coherent nontrivial modes are first order (Sec.~\ref{sec:modes}): $\Phi^{L}_1\propto h\,\xi_1$ and $\Phi^{R}_1\propto h\,\eta_1$, whose sign sets $\{\xi_1\gtrless0\}$ and $\{\eta_1\gtrless0\}$ are hyperplane partitions. The partitions are thus read off from the signs of the projections onto $v_1$ and $u_1$, the direct analogue of the Laplacian's leading eigenvector in spectral clustering~\cite{fiedler1973algebraic}. We tested the algorithm on both stochastic linear saddle systems and nonlinear bistable dynamics with saddle fixed points. The numerically estimated leading predictive vector aligns with the theoretical expectation as the sample size increases (Fig.~\ref{linear cosine size} and the right panels of Fig.~\ref{nonlinear system}).


\subsection{Practical regime}

Conditioning an escaping process has a cost: the weights concentrate on rare surviving trajectories, so the effective sample size decays roughly as $e^{-c\,\Tr(\Lcal)\,T}$. The estimator is thus reliable for moderate horizons and modest escape rates. The window-independence guaranteed by the Riccati attractor (Sec.~\ref{sec:riccati}) provides a plateau in $T$ over which the estimate is unbiased before the variance penalty dominates; one selects $T$ at the short end of the plateau.

\section{Special case: single, dominant unstable mode}
\label{sec:singlemode}

A case of particular physical importance is a saddle with a single unstable direction whose escape rate is much lower than the relaxation rate of every stable direction, $\lambda^+\ll\min_i|\lambda_i^-|$ (equivalently, $\lambda^+<\min_i|\lambda_i^-|$ and the long-lag
limit $\tau\to\infty$). This is precisely the separation of timescales underlying Kramers' and Eyring's classical
reaction-rate theory~\cite{hanggi1990reaction} and transition-path
theory~\cite{vanden2006towards}: a rare, slow escape along a single reaction coordinate against a fast-relaxing bath of transverse
modes. Conventionally, the reaction coordinate is \emph{identified} with the unstable direction of the linearized dynamics at the transition state~\cite{kramers1940brownian, langer1969statistical, hanggi1990reaction}; we show here that in this regime, the full saddle construction of
Secs.~\ref{sec:stab}--\ref{sec:datadriven} justifies that
identification with exact closed-form solutions.

The left and right singular functions of the Koopman operator correspond to predictive and retrospective coherent sets, respectively. Given the time-reversal symmetry, we focus on prediction in the rest of this section. Write $\ell^+,r^+$ for the unique unstable left and right eigenvectors of $A$, $\lambda^+$ for the unstable eigenvalue, and $\kappa\coloneq \ell^{+\top}\Sigma\ell^+$. The stabilizing solution (Eq.~\eqref{eq:hP}) is rank one, $P_+=\ell^+\ell^{+\top}/\kappa$.

\begin{proposition}[Left eigenvectors of $\Atil$ for single unstable mode]
\label{prop:singleleft}
For a saddle $A$ with a single unstable mode,
\begin{equation}
\begin{aligned}
&\ell^{+\top}\Atil=-\lambda^+\,\ell^{+\top},
\qquad
\widetilde\ell_i^{-\top}\Atil=\lambda_i^-\,\widetilde\ell_i^{-\top},\\
&\widetilde\ell_i^-\coloneq \ell_i^--\frac{\ell_i^{-\top}\Sigma\ell^+}{\kappa}\,\ell^+.
\end{aligned}
\label{eq:singleleft}
\end{equation}
In particular, the unstable left eigenvector of $A$ is a left eigenvector of $\Atil$ with flipped eigenvalue.
\end{proposition}
\noindent Proof in Appendix S10. Proposition~\ref{prop:reflect} already gives the matching statement for right eigenvectors ($\Atil r_i^-=\lambda_i^-r_i^-$, $\Atil\Sigma\ell^+=-\lambda^+\Sigma\ell^+$); Eq.~\eqref{eq:singleleft} is its left-eigenvector counterpart. The flip preserves $\ell^+$ \emph{exactly} whenever there
is only one unstable mode to flip. With the additional assumption of timescale separation such that the unstable mode dominates the long-lag limit of the stabilized transfer operator, the leading local predictive direction simply converges to $\ell^+$, which is orthogonal to the stable subspace (Fig.~\ref{linear demo}).

\begin{proposition}[Long-lag convergence of the leading predictive direction]
\label{prop:committor}
If $\lambda^+<\min_i|\lambda_i^-|$, then $-\lambda^+$ is the eigenvalue of
$\Atil$ closest to zero, so $\ell^+$ is the left eigenvector of the slowest
mode of $\Atil$ and
\begin{equation}
v_1(\tau)\ \xrightarrow[\ \tau\to\infty\ ]{}\ \ell^+
\qquad\text{(up to normalization).}
\label{eq:vlimit}
\end{equation}
\end{proposition}
Proof in Appendix S10. We now interpret the physical meaning of Proposition~\ref{prop:committor} by connecting it to the stable manifold and the committor. The stable manifold of a fixed point $\mathbf x_*$ is the set of points whose deterministic orbit converges to $\mathbf x_*$ rather than escaping~\cite{perko2013differential}. In the single-unstable regime of this section, the stable manifold partitions the phase space into two disjoint regions whose states possess distinct asymptotic future behaviors. For a linear saddle, the hyperplane orthogonal to the left unstable eigenvector, $\{\mathbf{x}:(\mathbf{x}-\mathbf{x}_*)^\top \ell^+=0\}$, is exactly the stable manifold. For a saddle fixed point in a genuinely nonlinear system, this hyperplane is tangent to the stable manifold, serving as a local linear approximation. Thus, the projection of states near the saddle onto the local predictive direction $v_1$ forecasts the asymptotic future evolution. We numerically verified that the estimated $\hat{v}_1$ is perpendicular to the stable manifold at the saddle (middle panels of Fig.~\ref{nonlinear system}).

For a stochastic system, the mathematical object that more directly quantifies the probabilistic future evolution is the committor $q(\mathbf{x})$, which is defined as the probability of reaching basin A before basin B conditioned on the initial state $\mathbf x$~\cite{vanden2006towards}. The committor depends on the global definitions of the two basins and therefore need not, in general, be determined by the local linearization. Nevertheless, for a linear saddle with basin boundaries compatible with the unstable coordinate, the committor is exactly an affinely normalized error function of $\ell^{+\top}(\mathbf{x}-\mathbf{x}^*)$, and its committor contours are hyperplanes normal to $\ell^+$ (Appendix S11 and Ref.~\cite{berezhkovskii2005one}). We empirically observed that the local predictive direction is orthogonal to the contour of the committor at the saddle (middle panels of Fig.~\ref{nonlinear system}). These results verify that the singular functions we construct, along with the numerical estimator, can yield local predictive representations consistent with classical theories.

The single-unstable case also closes the loop on the stabilized covariance with an inverse-free formula:

\begin{proposition}[Stabilized covariance for single unstable mode]
\label{prop:singlesigma}
For a saddle $A$ with a single unstable mode,
\begin{equation}
\Stil=\Sigma-\Sigma\ell^+r^{+\top}-r^+\ell^{+\top}\Sigma.
\label{eq:singlesigma}
\end{equation}
\end{proposition}
\noindent Proof in Appendix S10. Eq.~\eqref{eq:singlesigma} is the covariance-level flip of the unstable mode. As a check, in the all-unstable corner ($d=1$, $r^+=\ell^+=1$) this reduces to $\Stil=\Sigma-2\Sigma=-\Sigma$, matching Eq.~\eqref{eq:corner-unstable}. We note that Eqs.~\eqref{eq:singleleft} and \eqref{eq:singlesigma} do not extend to multiple unstable modes by the naive matrix analogue.

Equations~\eqref{eq:vlimit} and~\eqref{eq:singlesigma} together give $v_1,\Stil$ in closed form. In chemical reaction rates, this provides the stabilized transverse covariance directly, wherever transition-state theory's separation of timescales holds~\cite{kramers1940brownian,langer1969statistical,hanggi1990reaction}. In enhanced-sampling methods, a central practical problem is identifying the reaction coordinate to bias simulations along~\cite{laio2002escaping,valsson2016enhancing}, and our algorithm provides a data-driven estimator for that coordinate with a proven guarantee of what it converges to, rather than a heuristic collective variable chosen by hand~\cite{tiwary2016spectral}. The same closed form applies to any metastable system with this structure, independent of the specific physical mechanism.

\section{Discussion and connections}
\label{sec:discussion}

\emph{Transition-path theory.} A saddle is a transition state; its unstable manifold is the linearized reaction coordinate, and the conditioned process is the local linear theory of the reactive ensemble. When the saddle has a single unstable direction whose escape is the slowest process in the system---the classical rare-event regime of Kramers' and Eyring's rate theory~\cite{hanggi1990reaction,vanden2006towards}---the leading singular function is aligned with the linear committor direction (Sec.~\ref{sec:singlemode}).

\emph{Linear control and stochastic realization.} The stabilization solution is the attractor of an algebraic Riccati flow~\cite{lancaster1995algebraic}, and the eigenvalue reflection of Proposition~\ref{prop:reflect} is familiar from acausal stochastic realization~\cite{pavon1984canonical,akaike1975markovian}. The agreement is not only formal: $\Atil$ is the optimal closed-loop generator of the linear-quadratic regulator, and the same solution is selected there because the regulated dynamics must not diverge---the control-theoretic counterpart of requiring the conditioned process not to escape. Canonical correlations between past and future of a \emph{stable} stochastic system are likewise classical~\cite{akaike1975markovian,pavon1984canonical}; what is new here is the escaping case, where no invariant probability measure exists, and the correlation problem must first be repaired by conditioning.

\emph{Stochastic thermodynamics.} The selective time reversal of Proposition~\ref{prop:reflect} reverses the time-irreversible part of the drift only on the unstable manifold, flipping the entropy-producing arrow on the escaping directions while leaving the recurrent ones untouched---a subspace-resolved statement about dissipation.

\emph{Generative diffusion models.} The all-unstable case, Eq.~\eqref{eq:corner-unstable}, is exactly the reverse-time score-based SDE of generative sampling~\cite{anderson1982reverse,song2020score}. It is the extreme point of a family whose interior is the saddle, with conditioning acting only on the escaping subspace: reverse-time generative dynamics and transport through a transition state are two ends of one construction.

\emph{Escape rate and survival.} The decay $e^{-\Tr(\Lcal)\tau}$ sets the rate at which coherence is lost, and $(\Lb,h,-\Tr(\Lcal))$ is the principal eigen-triple of a killed (conditioned) semigroup: $h$ is the survival function, small where escape is imminent, and $\rho_{\Atil}=h\nu$ is the stationary law of the associated $Q$-process, with $\nu$ its dual eigenfunction living on the unstable side~\cite{collet2013quasi}. The same mechanism underlies the metastable spectra of open systems more broadly: a decay rate becomes visible only once an invariant probability measure is given up in favor of a survival-weighted or otherwise non-normalizable reference, as in the Pollicott--Ruelle resonances of chaotic flows and the Gamow states of quantum scattering~\cite{ruelle1986resonances,bohm1989dirac}. Here that role is played by the reference weight $\rhotil_0$ itself.

\emph{Limitations.} Real eigenvalues are assumed for clarity only; the complex case extends by flipping conjugate pairs together (Appendix S4 D). Timescale separation is needed for the committor correspondence, not for the closed forms. Linearity is a genuine restriction: the estimator applies to nonlinear systems only within a neighborhood whose radius must be established empirically, and its variance grows with the escape rate, limiting the usable horizons. With several unstable directions, the flip preserves the unstable left subspace but distinguishes no direction within it, leaving open what plays the role of the reaction coordinate. The framework is local by design, complementing global numerical coherent-set computations.

\section*{Acknowledgment}
We thank Joshua L. Pughe-Sanford, Charles L. Epstein, and Ruoyu Shang for helpful discussions and comments, which helped improve this paper. The Flatiron Institute is a division of the Simons Foundation.

\bibliography{refs}

\end{document}


\title{\texorpdfstring{Supplemental Material:\\ Local Koopman modes of saddle-point diffusions}{title}}
\author{Xuehao Ding}
\affiliation{Flatiron Institute, Simons Foundation}

\maketitle

We give full derivations for the results of the main text. Section~\ref{app:doob} derives the Doob-transformed generator; Sec.~\ref{app:gaussianeig} treats the Gaussian eigenfunction; Sec.~\ref{app:reflect} proves the eigenvalue reflection; Sec.~\ref{app:family} classifies the Riccati solutions; Sec.~\ref{app:flow} establishes the Riccati flow and the unique attractor; Sec.~\ref{app:Stil} proves the closed-form stabilized covariance; Sec.~\ref{app:unitary} proves the unitary equivalence; Sec.~\ref{app:singfun} proves the closed-form singular functions; Sec.~\ref{app:spaces} makes the function spaces rigorous; Sec.~\ref{app:single} proves the closed-form results for a single, dominant unstable mode; Sec.~\ref{app:linear_committor} relates the committor of a linear saddle to the unstable left eigenvector.

\section{The Doob transform and conditioning}
\label{app:doob}

We introduce the general (time-dependent) Doob transform, show that conditioning on a terminal distribution is a special case, and finally show that its long-horizon limit is the time-independent eigenfunction transform used in the main text.

\subsection{General \texorpdfstring{$h$}{h}-transform via a space-time harmonic function}

Let $h(\mathbf x,t)>0$ be a \emph{space-time harmonic function} of the diffusion Eq.~\eqref{main-eq:ou},
\begin{equation}
(\partial_t + \Lb)\,h(\mathbf x,t) = 0,
\label{eq:harmonic}
\end{equation}
equivalently $h(\mathbf x(t),t)$ is a martingale, $\E[h(\mathbf x(t'),t')\mid \mathbf x(t)=\mathbf x]=h(\mathbf x,t)$ for $t'>t$. Define the conditioned process by the tilted transition density
\begin{equation}
p^{h}(\mathbf y,t'\mid \mathbf x,t) = \frac{h(\mathbf y,t')}{h(\mathbf x,t)}\,p(\mathbf y,t'\mid \mathbf x,t).
\label{eq:tilted}
\end{equation}
This is a genuine transition density: it is positive, and the martingale property makes it normalized,
\begin{equation}
\int p^{h}(\mathbf y,t'\mid \mathbf x,t)\,d\mathbf y
= \frac{\E[h(\mathbf x(t'),t')\mid \mathbf x(t)=\mathbf x]}{h(\mathbf x,t)} = 1.
\end{equation}
Its generator follows by expanding Eq.~\eqref{eq:tilted} to first order in $\varepsilon=t'-t$. For a time-independent test function $g$, with $F(\mathbf x,t)\coloneq h(\mathbf x,t)g(\mathbf x)$,
\begin{equation}
\E^{h}\!\big[g(\mathbf x(t+\varepsilon))\mid \mathbf x(t)=\mathbf x\big]
= \frac{\E[F(\mathbf x(t+\varepsilon),t+\varepsilon)\mid \mathbf x(t)=\mathbf x]}{h(\mathbf x,t)}
= g + \frac{\varepsilon}{h}\,(\partial_t+\Lb)(hg) + o(\varepsilon),
\end{equation}
so the conditioned backward generator is 
\begin{equation}
\begin{aligned}
\Lb_{h}g&=\tfrac1h(\partial_t+\Lb)(hg)\\
&= \tfrac1h\Lb(hg)-\tfrac1h\Lb(h)g.
\end{aligned}
\end{equation}
Using the product rule $\Lb(hg)=g\,\Lb h + h\,\Lb g + D\nabla h\cdot\nabla g$ and the harmonic condition Eq.~\eqref{eq:harmonic},
\begin{equation}
\Lb_{h}g = \frac1h\big[\,g\,(\partial_t h+\Lb h) + h\,\Lb g + D\nabla h\cdot\nabla g\,\big]
= \Lb g + D\,\nabla\log h\cdot\nabla g .
\end{equation}
The diffusion is unchanged and the drift gains $D\nabla\log h(\mathbf x,t)$, which is Eq.~\eqref{main-eq:driftshift-point}. In general $h$ depends on $t$, so the added drift does too.

\subsection{Conditioning on a terminal distribution}

Fix a horizon $T$ and a target terminal density $\pi$. Let $p_T$ be the marginal density of the original process at time $T$, and set the reweighting $w=\pi/p_T$. Choose
\begin{equation}
h(\mathbf x,t) = \E\big[w(\mathbf x(T))\mid \mathbf x(t)=\mathbf x\big]
= \int p(\mathbf x_T,T\mid \mathbf x,t)\,\frac{\pi(\mathbf x_T)}{p_T(\mathbf x_T)}\,d\mathbf x_T
= (\Kop_{T-t}w)(\mathbf x),
\label{eq:hcond}
\end{equation}
which is Eq.~\eqref{main-eq:hfinal}. As a conditional expectation of a fixed terminal random variable, $h(\mathbf x(t),t)$ is automatically a martingale, so Eq.~\eqref{eq:harmonic} holds and the previous subsection applies.

This is exactly the transform that reshapes the terminal law into $\pi$. The path-measure reweighting is $dQ/dP|_{\mathcal F_t}=h(\mathbf x(t),t)/h(\mathbf x(0),0)$; at $t=T$, since $h(\mathbf x,T)=w(\mathbf x)$,
\begin{equation}
\frac{dQ}{dP}\Big|_{\mathcal F_T} \propto w(\mathbf x(T)) = \frac{\pi(\mathbf x(T))}{p_T(\mathbf x(T))},
\end{equation}
so the terminal marginal of $Q$ is $w(\mathbf x_T)\,p_T(\mathbf x_T)=\pi(\mathbf x_T)$. Conditioning the dynamics to end in the distribution $\pi$ is therefore precisely the Doob transform with $h$ given by Eq.~\eqref{eq:hcond}.

\subsection{Long-horizon limit: the autonomous eigenfunction transform}

For a Gaussian weight $w(\mathbf x)\propto\exp(-\tfrac12\mathbf x^{\!\top}\Sigma'^{-1}\mathbf x)$, the OU semigroup keeps $h_s\coloneq\Kop_s w$ Gaussian, $h_s(\mathbf x)\propto C_s\exp(\tfrac12\mathbf x^{\!\top}P_s\mathbf x)$ with $s=T-t$. Substituting into $\partial_s h_s=\Lb h_s$ separates into the Riccati flow $\dot P_s=A^{\!\top}P_s+P_sA+P_sDP_s$ (Sec.~\ref{app:flow}) for the quadratic part and $\tfrac{d}{ds}\log C_s=\tfrac12\Tr(DP_s)$ for the scalar. As $s\to\infty$, $P_s\to P_+=\Lp(\Lp^{\!\top}\Sigma\Lp)^{-1}\Lp^{\!\top}$, the stabilizing attractor (Sec.~\ref{app:flow}), so the conditioning drift becomes \emph{autonomous},
\begin{equation}
D\nabla\log h(\mathbf x,t)=DP_{T-t}\,\mathbf x \;\xrightarrow[\;T-t\to\infty\;]{}\; DP_+\mathbf x,
\qquad \Atil=A+DP_+.
\end{equation}
The limiting spatial profile $h(\mathbf x)=\exp(\tfrac12\mathbf x^{\!\top}P_+\mathbf x)$ is a time-independent eigenfunction, $\Lb h=\lambda h$ with $\lambda=\tfrac12\Tr(DP_+)=-\Tr(\Lcal)$ (Sec.~\ref{app:gaussianeig}); the associated space-time harmonic function is $h(\mathbf x,t)=e^{-\lambda t}h(\mathbf x)$, which indeed satisfies Eq.~\eqref{eq:harmonic}. Because only $\nabla\log h(\mathbf x )$ enters the dynamics, the scalar $C_s$ (equivalently $e^{-\lambda t}$) is dynamically irrelevant, and the conditioned process is the time-homogeneous ($\lambda$-invariant) Doob transform, with transition kernel and intertwining relation
\begin{equation}
p^{h}_\tau(\mathbf x,\mathbf y)=e^{-\lambda\tau}\,\frac{h(\mathbf y)}{h(\mathbf x)}\,p_\tau(\mathbf y\mid\mathbf x),
\qquad
\Kop_\tau=e^{\lambda\tau}\,M_h\,\widetilde\Kop_\tau\,M_h^{-1},\quad \lambda=-\Tr(\Lcal),
\label{eq:autonomous}
\end{equation}
where $M_h$ is multiplication by $h$ and $\widetilde\Kop_\tau$ is the Koopman operator of the stabilized process. This time-homogeneous transform is what the main text and the remaining appendices use; its legitimacy as a true martingale requires the finite exponential moment $\E_{\mathbf x}[h(\mathbf x(\tau))]<\infty$, discussed in Sec.~\ref{app:spaces}.

\section{The Gaussian eigenfunction and its eigenvalue}
\label{app:gaussianeig}

For $h=\exp(\tfrac12 x^{\!\top}Px)$ ($P=P^{\!\top}$), with $\nabla h=(Px)h$ and $\partial_\alpha\partial_\beta h=(P_{\alpha\beta}+(Px)_\alpha(Px)_\beta)h$,
\begin{equation}
\Lb h = \Big[\tfrac12 x^{\!\top}\!\big(A^{\!\top}P+PA+PDP\big)x+\tfrac12\Tr(DP)\Big]h .
\end{equation}
The quadratic term vanishes iff $P$ solves the algebraic Riccati equation
\begin{equation}
A^{\!\top}P+PA+PDP=0,
\label{eq:riccatiS}
\end{equation}
and then $\Lb h=\lambda h$ with $\lambda=\tfrac12\Tr(DP)$. For $P=\Lp(\Lp^{\!\top}\Sigma\Lp)^{-1}\Lp^{\!\top}$, put $S=\Lp^{\!\top}\Sigma\Lp$; from Eq.~\eqref{main-eq:lyap} and $A^{\!\top}\Lp=\Lp\Lcal$,
\begin{equation}
\Lp^{\!\top}D\Lp=-\Lp^{\!\top}(A\Sigma+\Sigma A^{\!\top})\Lp=-\Lcal S-S\Lcal,
\qquad
\lambda=\tfrac12\Tr\!\big[(-\Lcal S-S\Lcal)S^{-1}\big]=-\Tr(\Lcal).
\end{equation}
Hence $\Kop_\tau h=e^{-\Tr(\Lcal)\tau}h$.

\section{Proof of the eigenvalue reflection (main-text Proposition~1)}
\label{app:reflect}

Let $\Atil=A+DP_+$, $P_+=\Lp(\Lp^{\!\top}\Sigma\Lp)^{-1}\Lp^{\!\top}$. For a stable right eigenvector, biorthogonality gives $\Lp^{\!\top}r_i^-=0$, so $P_+r_i^-=0$ and $\Atil r_i^-=\lambda_i^- r_i^-$. For an unstable left eigenvector $\ell_i^+$, $\Lp^{\!\top}\Sigma\ell_i^+=Se_i$ gives $P_+\Sigma\ell_i^+=\ell_i^+$, whence
\begin{equation}
\Atil(\Sigma\ell_i^+)=A\Sigma\ell_i^+ + D\ell_i^+
= A\Sigma\ell_i^+ - (A\Sigma+\Sigma A^{\!\top})\ell_i^+
= -\Sigma A^{\!\top}\ell_i^+ = -\lambda_i^+(\Sigma\ell_i^+),
\end{equation}
using Eq.~\eqref{main-eq:lyap}. Thus $\mathrm{spec}(\Atil)=\{\lambda_i^-\}\cup\{-\lambda_i^+\}$, all negative, so $\Atil$ is Hurwitz. $\qquad\blacksquare$

\section{The Riccati solution family and the Hamiltonian picture}
\label{app:family}

We derive the solution family (Eq.~\eqref{main-eq:Pfamily}) of the algebraic Riccati equation
\begin{equation}
A^{\!\top}P+PA+PDP=0
\label{eq:riccatiS4}
\end{equation}
and identify the stabilizing solution, first for real spectra and then for complex ones. Since $P$ enters through the quadratic form, we restrict throughout to symmetric solutions $P=P^\top$.

\subsection{Riccati solutions as invariant subspaces}

$P$ solves Eq.~\eqref{eq:riccatiS4} iff the \emph{graph subspace} $\mathcal V_P\coloneq\mathrm{Im}\bigl[\begin{smallmatrix}I\\P\end{smallmatrix}\bigr] =\{(\mathbf z,P\mathbf z):\mathbf z\in\R^d\}\subset\R^{2d}$ is invariant under the Hamiltonian matrix $\mathcal H=\bigl[\begin{smallmatrix}A&D\\0&-A^{\!\top}\end{smallmatrix}\bigr]$, in which case $\mathcal H$ acts on $\mathcal V_P$ as $\Atil=A+DP$. Indeed, $\mathcal H\bigl[\begin{smallmatrix}I\\P\end{smallmatrix}\bigr] =\bigl[\begin{smallmatrix}A+DP\\-A^{\!\top}P\end{smallmatrix}\bigr]$, whose bottom block equals $P$ times its top block iff $-A^{\!\top}P=P(A+DP)$, i.e.\ Eq.~\eqref{eq:riccatiS4}; the top block $A+DP$ is then the action of $\mathcal H$ on $\mathcal V_P$ in the coordinate $\mathbf z$.

\subsection{Eigenvectors of \texorpdfstring{$\mathcal H$}{H} (real spectrum)}

Being block upper-triangular, $\mathcal H$ has the $2d$ eigenvalues of $A$ together with those of $-A^{\!\top}$, i.e.\ the mirror pairs $\{\lambda_i,-\lambda_i\}$, with eigenvectors
\begin{equation}
\mathcal H\begin{bmatrix}r_i\\0\end{bmatrix}=\lambda_i\begin{bmatrix}r_i\\0\end{bmatrix},
\qquad
\mathcal H\begin{bmatrix}\Sigma\ell_i\\ \ell_i\end{bmatrix}=-\lambda_i\begin{bmatrix}\Sigma\ell_i\\ \ell_i\end{bmatrix}.
\label{eq:Heig}
\end{equation}
The first is immediate. For the second, the bottom block gives $-A^{\!\top}\ell_i=-\lambda_i\ell_i$,
and the top block, using $D=-(A\Sigma+\Sigma A^{\!\top})$ and $A^{\!\top}\ell_i=\lambda_i\ell_i$, is
\begin{equation}
A\Sigma\ell_i+D\ell_i=A\Sigma\ell_i-(A\Sigma+\Sigma A^{\!\top})\ell_i=-\Sigma A^{\!\top}\ell_i=-\lambda_i\,\Sigma\ell_i,
\end{equation}
confirming eigenvalue $-\lambda_i$ with eigenvector $[\Sigma\ell_i;\ell_i]$.

\subsection{Graph subspaces: keep or flip}

For a symmetric Riccati solution, the invariant graph $\mathcal V_P$ is Lagrangian with respect to the standard symplectic form defined by
$
J=
\begin{bmatrix}
0&I\\
-I&0
\end{bmatrix}.
$
The two Hamiltonian eigenvectors belonging to the same mirror pair,
$z_i^+=
\begin{bmatrix}
r_i\\
0
\end{bmatrix}\text{ and }
z_i^- =
\begin{bmatrix}
\Sigma \ell_i\\
\ell_i
\end{bmatrix},
$
have nonzero symplectic product,
$
(z_i^+)^\top J z_i^- = r_i^\top \ell_i = 1,
$
and therefore cannot both belong to the same Lagrangian subspace. Since $\mathcal V_P$ has dimension $d$, it must contain exactly one
eigenvector from each mirror pair $\{\lambda_i,-\lambda_i\}$. Let $\mathcal F\subseteq\{1,\dots,d\}$ index the pairs for which we take the $-\lambda_i$ eigenvector (``flip''); the complement takes the $+\lambda_i$ eigenvector (``keep''). The spanning set is
\begin{equation}
\big\{\,[r_i;0]:i\notin\mathcal F\,\big\}\ \cup\ \big\{\,[\Sigma\ell_i;\ell_i]:i\in\mathcal F\,\big\}.
\end{equation}

The matrix $P$ representing this subspace as a graph is fixed by the defining relation of $\mathcal V_P$: every spanning vector $[\mathbf a;\mathbf b]$ must satisfy $\mathbf b=P\mathbf a$. Reading off the two vector types gives the $d$ conditions
\begin{equation}
P\,r_i = 0\quad(i\notin\mathcal F),
\qquad
P\,(\Sigma\ell_i) = \ell_i\quad(i\in\mathcal F).
\label{eq:Pconditions}
\end{equation}
The top parts $\{r_i\}_{i\notin\mathcal F}\cup\{\Sigma\ell_i\}_{i\in\mathcal F}$ form a basis of $\R^d$, so Eq.~\eqref{eq:Pconditions} determines $P$ uniquely. The solution is
\begin{equation}
P_{\mathcal F}=L_{\mathcal F}\big(L_{\mathcal F}^{\!\top}\Sigma L_{\mathcal F}\big)^{-1}L_{\mathcal F}^{\!\top},
\qquad L_{\mathcal F}=[\ell_i]_{i\in\mathcal F},
\label{eq:PfamilyS4}
\end{equation}
which we verify against Eq.~\eqref{eq:Pconditions}. For $i\notin\mathcal F$, biorthogonality $\ell_j^{\!\top}r_i=\delta_{ji}=0$ (all $j\in\mathcal F$) gives $L_{\mathcal F}^{\!\top}r_i=0$, hence $P_{\mathcal F}r_i=0$. For $i\in\mathcal F$, the vector $L_{\mathcal F}^{\!\top}\Sigma\ell_i$ has entries $[\ell_j^{\!\top}\Sigma\ell_i]_{j\in\mathcal F}=(L_{\mathcal F}^{\!\top}\Sigma L_{\mathcal F})\mathbf e_i$, so
\begin{equation}
P_{\mathcal F}\,\Sigma\ell_i
= L_{\mathcal F}\big(L_{\mathcal F}^{\!\top}\Sigma L_{\mathcal F}\big)^{-1}
(L_{\mathcal F}^{\!\top}\Sigma L_{\mathcal F})\mathbf e_i
= L_{\mathcal F}\mathbf e_i=\ell_i,
\end{equation}
as required. The action $\Atil_{\mathcal F}=A+DP_{\mathcal F}$ then has eigenvalues $\{\lambda_i:i\notin\mathcal F\}\cup\{-\lambda_i:i\in\mathcal F\}$: kept eigenvalues are unchanged, flipped ones negated. Choosing $\mathcal F=\{\text{unstable modes}\}$ negates exactly the positive eigenvalues, making $\Atil$ Hurwitz---the stabilizing solution, Eq.~\eqref{main-eq:hP}.

\subsection{Complex spectrum}
\label{app:complex}

If $A$ has a complex eigenvalue $\lambda=a+ib$ ($b\neq0$), then $\bar\lambda=a-ib$ is also an eigenvalue, with conjugate eigenvectors $r_{\bar\lambda}=\overline{r_\lambda}$, $\ell_{\bar\lambda}=\overline{\ell_\lambda}$. The Hamiltonian then carries the four eigenvalues $\{\lambda,\bar\lambda,-\lambda,-\bar\lambda\}$: two mirror pairs exchanged by conjugation. For $P_{\mathcal F}$---and hence $\Atil$ and the reference measure---to be \emph{real}, the flipped set $\mathcal F$ must be closed under conjugation: a complex mode is flipped only together with its conjugate. One therefore flips \emph{spiral modes} (conjugate pairs) as a unit, never an individual complex eigenvalue.

Concretely, for a flipped conjugate pair the relevant $\mathcal H$-invariant subspace is the real
plane
\begin{equation}
\mathrm{span}_{\R}\Big\{\ \mathrm{Re}\begin{bmatrix}\Sigma\ell_\lambda\\ \ell_\lambda\end{bmatrix},\
\mathrm{Im}\begin{bmatrix}\Sigma\ell_\lambda\\ \ell_\lambda\end{bmatrix}\ \Big\},
\end{equation}
and $P_{\mathcal F}$ from Eq.~\eqref{eq:PfamilyS4} with $L_{\mathcal F}$ containing the pair $\{\ell_\lambda,\bar\ell_\lambda\}$ is real, because the conjugation-symmetric combination $L_{\mathcal F}(L_{\mathcal F}^{\!\top}\Sigma L_{\mathcal F})^{-1}L_{\mathcal F}^{\!\top}$ has canceling imaginary parts. Under the flip, the eigenvalue pair $\{a\pm ib\}$ of $A$ becomes $\{-a\pm ib\}$ in $\Atil$: the growth rate $a$ is reflected to $-a$, while the oscillation frequency $b$ is preserved. An \emph{unstable spiral} ($a>0$) thus becomes a \emph{stable spiral} ($-a<0$) of the same frequency.
\section{Riccati flow and uniqueness of the attractor}
\label{app:flow}

From Eq.~\eqref{main-eq:hfinal}, $h_s=\Kop_s w$ with $w\propto\exp(-\tfrac12 x^{\!\top}\Sigma'^{-1}x)$. The Gaussian class is preserved, so $h_s\propto\exp(\tfrac12 x^{\!\top}P_s x)$, and $\partial_s h_s=\Lb h_s$ with the quadratic-form matching of Sec.~\ref{app:gaussianeig} gives
\begin{equation}
\dot P_s = A^{\!\top}P_s+P_sA+P_sDP_s = \mathcal R(P_s).
\label{eq:flowS}
\end{equation}

We prove that $P_+=\Lp(\Lp^{\!\top}\Sigma\Lp)^{-1}\Lp^{\!\top}$ is the unique attractor of Eq.~\eqref{eq:flowS} on the symmetric matrix subspace; every other solution is a saddle.
\begin{proof}
With $E=P-P_{\mathcal F}$ and $\mathcal R(P_{\mathcal F})=0$,
\begin{equation}
\dot E = A^{\!\top}E+EA+P_{\mathcal F}DE+EDP_{\mathcal F}+EDE
= \Atil_{\mathcal F}^{\!\top}E+E\Atil_{\mathcal F}+EDE,\qquad \Atil_{\mathcal F}=A+DP_{\mathcal F}.
\label{eq:deviation}
\end{equation}
The linearization $E\mapsto\Atil_{\mathcal F}^{\!\top}E+E\Atil_{\mathcal F}$ as a Sylvester matrix ODE has eigenvalues $\widetilde\lambda_i+\widetilde\lambda_j$, where $\{\widetilde\lambda_i\}=\mathrm{spec}(\Atil_{\mathcal F})$. For the stabilizing solution all $\widetilde\lambda_i<0$, so all $\widetilde\lambda_i+\widetilde\lambda_j<0$ and $P_+$ is locally exponentially stable; for any other $\mathcal F$ some $\widetilde\lambda_i>0$, giving an unstable direction. Globally, with $\Phi_s=e^{\Atil s}$ ($\Atil$ Hurwitz), Eq.~\eqref{eq:deviation} integrates to
\begin{equation}
E_s=\Phi_s^{\!\top}\Big(E_0^{-1}-\int_0^s\Phi_uD\Phi_u^{\!\top}\,du\Big)^{-1}\Phi_s\xrightarrow{s\to\infty}0,
\end{equation}
for any admissible $E_0$, independently of $\Sigma'$.
\end{proof}

\section{Closed form for the stabilized covariance at an arbitrary flip set}
\label{app:Stil}
We prove the following more general version of Proposition~\ref{main-prop:Stil}: for any flip set $\mathcal F\subseteq\{1,\dots,d\}$ with complement $\mathcal F^c$,
\begin{equation}
\Stil_{\mathcal F}=(P_{\mathcal F^c}-P_{\mathcal F})^{-1},
\end{equation}
where $\Stil_{\mathcal F}$ is the solution of the Lyapunov equation $\Atil_{\mathcal F}\Stil_{\mathcal F}+\Stil_{\mathcal F}\Atil_{\mathcal F}^{\!\top}=-D$ with $\Atil_{\mathcal F}=A+DP_{\mathcal F}$.

\begin{proof}
Both $P_{\mathcal F}$ and $P_{\mathcal F^c}$ solve the Riccati equation $A^{\!\top}P+PA+PDP=0$. Set $M\coloneq P_{\mathcal F^c}-P_{\mathcal F}$ and expand, using $M+P_{\mathcal F}=P_{\mathcal F^c}$:
\begin{equation}
\begin{aligned}
M\Atil_{\mathcal F}+\Atil_{\mathcal F}^{\!\top}M+MDM
=&MA+A^{\!\top}M+\underbrace{MDP_{\mathcal F}+P_{\mathcal F}DM+MDM}_{=\,P_{\mathcal
F^c}DP_{\mathcal F^c}-P_{\mathcal F}DP_{\mathcal F}}\\
=&\underbrace{\big[A^{\!\top}P_{\mathcal F^c}+P_{\mathcal
F^c}A+P_{\mathcal F^c}DP_{\mathcal F^c}\big]}_{=0}
-\underbrace{\big[A^{\!\top}P_{\mathcal F}+P_{\mathcal
F}A+P_{\mathcal F}DP_{\mathcal F}\big]}_{=0}\\
=&0,
\end{aligned}
\label{eq:shift coordinate}
\end{equation}
so $M$ solves the same Riccati equation, based at $\Atil_{\mathcal F}$. So far, we have not used the complementary property. Eq.~\eqref{eq:shift coordinate} verifies that re-basing the Riccati equation at $\Atil_{\mathcal F}$ is a change of origin in the Riccati solution set, translating each solution $P_{\mathcal G}$ to $P_{\mathcal G}- P_{\mathcal F}$ with flip set $\mathcal F \triangle \mathcal G$ (symmetric difference). Since the rank of a Riccati solution is equal to the number of modes in the flip set, $P_{\mathcal G}- P_{\mathcal F}$ is full rank iff $\mathcal G = \mathcal F^c$. Thus, $M$ is invertible. Multiplying $M^{-1}$ to both sides of Eq.~\eqref{eq:shift coordinate} yields the Lyapunov equation:
\begin{equation}
    \Atil_{\mathcal F}M^{-1}+M^{-1}\Atil_{\mathcal F}^{\!\top}=-D.
\end{equation}
By the uniqueness of the solution, $\Stil_{\mathcal F}=M^{-1}$.
\end{proof}

\section{Proof of the unitary equivalence (main-text Proposition~\ref{main-prop:unitary})}
\label{app:unitary}

\emph{Isometry.} For $f\in L^2(\rho_{\Atil})$, $\|M_hf\|^2_{L^2(\rhotil_0)}=\int h^2f^2\rhotil_0\,dx=\int f^2\rho_{\Atil}\,dx=\|f\|^2_{L^2(\rho_{\Atil})}$, using $\rhotil_0=\rho_{\Atil}h^{-2}$; invertibility of $M_h$ ($h>0$) makes it unitary.

\emph{Intertwining.} From Sec.~\ref{app:doob}, $(\Kop_\tau(hg))(x)=\int p_\tau(y\mid x)h(y)g(y)\,dy=e^{\lambda\tau}h(x)(\widetilde\Kop_\tau g)(x)$, i.e.\ $\Kop_\tau M_h=e^{\lambda\tau}M_h\widetilde\Kop_\tau$, so $\Kop_\tau=e^{-\Tr(\Lcal)\tau}M_h\widetilde\Kop_\tau M_h^{-1}$.

\emph{Coincidence of singular systems.} If $\widetilde\Kop_\tau g_{\bm n}=\gamma^{\bm n}f_{\bm n}$ with $\{f_{\bm n}\},\{g_{\bm n}\}$ orthonormal in $L^2(\rho_{\Atil})$, then $\Kop_\tau(M_hg_{\bm n})=e^{-\Tr(\Lcal)\tau}\gamma^{\bm n}(M_hf_{\bm n})$ with $\{M_hf_{\bm n}\},\{M_hg_{\bm n}\}$ orthonormal in $L^2(\rhotil_0)$; hence the singular functions are $h f_{\bm n}, h g_{\bm n}$ and the singular values $e^{-\Tr(\Lcal)\tau}\gamma^{\bm n}$. The stabilized system $\{f_{\bm n},g_{\bm n}\}$ is what CCA computes (Sec.~\ref{main-sec:cca} of the main text): the CCA operator $C_0^{-1}C_\tau C_0^{-1}C_{-\tau}$ for the stabilized process equals $F=e^{\Atil^{\!\top}\tau}\Stil^{-1}e^{\Atil\tau}\Stil$. $\qquad\blacksquare$

\section{Closed-form singular functions}
\label{app:singfun}

We prove Eq.~\eqref{main-eq:mainresult}. By Sec.~\ref{app:unitary} it suffices to find the singular functions of the \emph{stabilized} process $\dot x=\Atil x+\xi$ ($\Atil$ Hurwitz), with stationary law $\rho=\mathcal N(0,\Stil)$, $\Atil\Stil+\Stil\Atil^{\!\top}=-D$; the saddle result then follows by multiplying by $h$ and rescaling by $e^{-\Tr(\Lcal)\tau}$. Throughout, $\He_n$ is the probabilists' Hermite polynomial, $\sum_n\frac{t^n}{n!}\He_n(z)=e^{tz-t^2/2}$, $\int\He_m\He_n\,d\mathcal N(0,1)=n!\,\delta_{mn}$, and $\bm n!=\prod_i n_i!$, $\gamma^{\bm n}=\prod_i\gamma_i^{n_i}$.

\begin{lemma}[univariate Mehler relation]
\label{lem:mehler}
If $(\xi,\eta)$ is standard bivariate normal with unit variances and $\E[\xi\eta]=\gamma$, then $\E[\He_n(\eta)\mid\xi]=\gamma^n\He_n(\xi)$.
\end{lemma}
\begin{proof}
Conditionally, $\eta\mid\xi\sim\mathcal N(\gamma\xi,1-\gamma^2)$, so
\begin{equation}
\E\!\big[e^{t\eta-t^2/2}\,\big|\,\xi\big]
=\exp\!\Big(t\gamma\xi+\tfrac12 t^2(1-\gamma^2)-\tfrac12 t^2\Big)
=\exp\!\Big((t\gamma)\xi-\tfrac12(t\gamma)^2\Big)
=\sum_{n\ge0}\frac{(t\gamma)^n}{n!}\He_n(\xi).
\end{equation}
Matching powers of $t$ with $\E[e^{t\eta-t^2/2}\mid\xi]=\sum_n\frac{t^n}{n!}\E[\He_n(\eta)\mid\xi]$ gives the claim.
\end{proof}

\paragraph{Step 1: canonical coordinates.} Under $\rho$, $(x_0,x_\tau)$ is jointly Gaussian, mean zero, with $\E[x_0x_0^{\!\top}]=\E[x_\tau x_\tau^{\!\top}]=\Stil$ and $\E[x_0x_\tau^{\!\top}]=\Stil e^{\Atil^{\!\top}\tau}$. Whitening $z_0=\Stil^{-1/2}x_0$, $z_\tau=\Stil^{-1/2}x_\tau$ (each $\mathcal N(0,I)$), the cross-covariance is $K=\Stil^{1/2}e^{\Atil^{\!\top}\tau}\Stil^{-1/2}$, real. Its SVD $K=\Psi\Gamma\Phi^{\!\top}$ ($\Gamma=\mathrm{diag}(\gamma_i)$, $\Psi,\Phi$ orthogonal) defines the canonical directions and coordinates
\begin{equation}
v_i=\Stil^{-1/2}\psi_i,\quad u_i=\Stil^{-1/2}\phi_i,\qquad
\xi_i=v_i^{\!\top}x_0,\quad \eta_i=u_i^{\!\top}x_\tau,
\end{equation}
for which
\begin{equation}
\E[\xi_i\xi_j]=\delta_{ij},\qquad \E[\eta_i\eta_j]=\delta_{ij},\qquad \E[\xi_i\eta_j]=\gamma_i\,\delta_{ij}.
\label{eq:canon}
\end{equation}
Since $Q_\tau=\Stil-e^{\Atil\tau}\Stil e^{\Atil^{\!\top}\tau}\succ0$ gives $K^{\!\top}K\prec I$, and $e^{\Atil\tau}$ is invertible, $\gamma_i\in(0,1)$. Note $K$ is real regardless of whether $\Atil$'s eigenvalues are real, so the construction holds for spiral saddles.

\paragraph{Step 2: independence and conditional decoupling.} Eqs.~\eqref{eq:canon} say the covariance of $(\xi_1,\eta_1,\dots,\xi_d,\eta_d)$ is block-diagonal in the pairs $(\xi_i,\eta_i)$, each block $\bigl[\begin{smallmatrix}1&\gamma_i\\\gamma_i&1\end{smallmatrix}\bigr]$, with zero covariance across pairs. Jointly Gaussian variables with zero covariance are independent, so the pairs $\{(\xi_i,\eta_i)\}$ are mutually independent. As $x_0\mapsto(\xi_1,\dots,\xi_d)$ is a linear bijection, conditioning on $x_0$ equals conditioning on $\xi$, and
\begin{equation}
\eta_i\mid\xi \ \sim\ \mathcal N(\gamma_i\xi_i,\,1-\gamma_i^2),\quad\text{with the }\eta_i\text{ conditionally independent};
\label{eq:conddecouple}
\end{equation}
in particular $\eta_i\mid\xi$ depends only on $\xi_i$.

\paragraph{Step 3: action on Hermite products.} Define $\mathcal H^v_{\bm n}(x)=\prod_i\He_{n_i}(v_i^{\!\top}x)$ and $\mathcal H^u_{\bm n}(x)=\prod_i\He_{n_i}(u_i^{\!\top}x)$. Then, using Eq.~\eqref{eq:conddecouple} and Lemma~\ref{lem:mehler} on each independent pair,
\begin{equation}
(\widetilde\Kop_\tau\mathcal H^u_{\bm n})(x)
=\E\!\Big[\prod_i\He_{n_i}(\eta_i)\,\Big|\,\xi\Big]
=\prod_i\E[\He_{n_i}(\eta_i)\mid\xi_i]
=\prod_i\gamma_i^{n_i}\He_{n_i}(\xi_i)
=\gamma^{\bm n}\,\mathcal H^v_{\bm n}(x).
\label{eq:exactaction}
\end{equation}

\paragraph{Step 4: orthonormality, completeness, and the singular system.} Under $\rho$ the $\xi_i=v_i^{\!\top}x$ are i.i.d.\ standard Gaussian ($v_i^{\!\top}\Stil v_j=\psi_i^{\!\top}\psi_j=\delta_{ij}$), so $\langle\mathcal H^v_{\bm m},\mathcal H^v_{\bm n}\rangle_\rho=\bm n!\,\delta_{\bm m\bm n}$, and likewise for $\mathcal H^u$; both families are orthogonal. They are complete because the univariate Hermite functions are a basis of $L^2(\mathcal N(0,1))$, so their products are a basis of $L^2(\mathcal N(0,I_d))$, pulled back to $L^2(\rho)$ by the bijection $x\mapsto\xi$. Writing $f_{\bm n}=\mathcal H^v_{\bm n}/\sqrt{\bm n!}$, $g_{\bm n}=\mathcal H^u_{\bm n}/\sqrt{\bm n!}$, Eq.~\eqref{eq:exactaction} reads $\widetilde\Kop_\tau g_{\bm n}=\gamma^{\bm n}f_{\bm n}$; by the $\xi\leftrightarrow\eta$ symmetry of Eq.~\eqref{eq:canon}, the adjoint satisfies $\widetilde\Kop_\tau^{\dagger}f_{\bm n}=\gamma^{\bm n}g_{\bm n}$. With orthonormality and completeness, this is by definition the singular value decomposition of the stabilized operator:
\begin{equation}
f_{\bm n}(x)=\frac{1}{\sqrt{\bm n!}}\prod_i\He_{n_i}(v_i^{\!\top}x),\quad
g_{\bm n}(x)=\frac{1}{\sqrt{\bm n!}}\prod_i\He_{n_i}(u_i^{\!\top}x),\quad
\text{singular values }\ \gamma^{\bm n}.
\end{equation}

\paragraph{Step 5: transport to the saddle case.} By Proposition~\ref{main-prop:unitary} (Sec.~\ref{app:unitary}), $\Kop_\tau=e^{-\Tr(\Lcal)\tau}M_h\widetilde\Kop_\tau M_h^{-1}$ with $M_h$ unitary, so applying $M_h$ and rescaling,
\begin{equation}
\Phi^L_{\bm n}=\frac{h}{\sqrt{\bm n!}}\prod_i\He_{n_i}(v_i^{\!\top}x),\quad
\Phi^R_{\bm n}=\frac{h}{\sqrt{\bm n!}}\prod_i\He_{n_i}(u_i^{\!\top}x),\quad
\sigma_{\bm n}=e^{-\Tr(\Lcal)\tau}\gamma^{\bm n},
\end{equation}
orthonormal in $L^2(\rhotil_0)$. This is Eq.~\eqref{main-eq:mainresult}. $\qquad\blacksquare$

\section{Function spaces and well-posedness}
\label{app:spaces}

\emph{Weighted $L^2$.} For a positive measurable weight $w$, $L^2(w)=\{f:\int|f|^2w\,dx<\infty\}$ is a Hilbert space whether or not $\int w\,dx<\infty$ (Lebesgue measure is the standard infinite example). Thus $\rhotil_0$ need not be normalizable to define $L^2(\rhotil_0)$; only the functions used must be square-integrable, and for Eq.~\eqref{main-eq:mainresult},
\begin{equation}
\|\Phi_{\bm n}\|^2_{L^2(\rhotil_0)}=\int h^2|\mathcal H_{\bm n}|^2\,\rhotil_0\,dx=\int|\mathcal H_{\bm n}|^2\,\rho_{\Atil}\,dx<\infty,
\end{equation}
using $\rhotil_0=\rho_{\Atil}h^{-2}$.

\emph{$\Kop_\tau$ is bounded on $L^2(\rhotil_0)$.} On $L^2(\rho_{\Atil})$, $\widetilde\Kop_\tau$ is Hilbert--Schmidt (Mehler-type) with discrete singular system $\sigma_{\bm n}\le1$. By the unitary equivalence of Proposition~\ref{main-prop:unitary}, $\Kop_\tau=e^{-\Tr(\Lcal)\tau}M_h\widetilde\Kop_\tau M_h^{-1}$ with $M_h:L^2(\rho_{\Atil})\to L^2(\rhotil_0)$ unitary, so $\Kop_\tau$ is itself bounded on $L^2(\rhotil_0)$, with operator norm $e^{-\Tr(\Lcal)\tau}<1$: a genuine contraction, not merely a formally defined object. 

\emph{True-martingale condition.} The change of measure of Sec.~\ref{app:doob} is legitimate iff $\E_x[h(X_\tau)]<\infty$, i.e.
\begin{equation}
\int e^{\frac12 y^{\!\top}Py}\,\mathcal N(y;e^{A\tau}x,Q_\tau)\,dy<\infty \iff P\prec Q_\tau^{-1},
\qquad Q_\tau=\int_0^\tau e^{As}De^{A^{\!\top}s}\,ds.
\end{equation}
This is automatic when $h$ is bounded/decaying; otherwise, it imposes a restriction on $\tau$. Since $\rho_{\Atil}$ is a genuine Gaussian ($\Atil$ Hurwitz, Sec.~\ref{app:reflect}), the problem on $L^2(\rho_{\Atil})$ is unconditionally well-posed. Because $h\le1$ (Sec.~\ref{main-sec:riccati}), this condition holds for every $\tau$ here; no restriction on the horizon is needed.

\section{Closed-form results for a single, dominant unstable mode}
\label{app:single}

Throughout this appendix, there is a single unstable direction: $\Lp=\ell^+$ (a vector, not a matrix), $\Lcal=\lambda^+$ (a scalar), $P_+=\ell^+\ell^{+\top}/\kappa$, $\kappa\coloneq\ell^{+\top}\Sigma\ell^+$.

\begin{lemma}
\label{lem:bilinear}
$\ell^{+\top}D\ell^+=-2\lambda^+\kappa$, and for any stable left eigenvector
$\ell_i^-$, $\ell_i^{-\top}D\ell^+=-(\lambda_i^-+\lambda^+)\,\ell_i^{-\top}\Sigma\ell^+$.
\end{lemma}
\begin{proof}
From $D=-(A\Sigma+\Sigma A^{\!\top})$ and $\ell^{+\top}A=\lambda^+\ell^{+\top}$, $A^{\!\top}\ell^+=\lambda^+\ell^+$:
\begin{equation}
\ell^{+\top}D\ell^+=-\ell^{+\top}A\Sigma\ell^+-\ell^{+\top}\Sigma A^{\!\top}\ell^+
=-\lambda^+\kappa-\lambda^+\kappa=-2\lambda^+\kappa.
\end{equation}
Replacing one $\ell^+$ by $\ell_i^-$ (using $\ell_i^{-\top}A=\lambda_i^-\ell_i^{-\top}$ on that factor) gives the second identity by the same substitution.
\end{proof}

\begin{proof}[Proof of Proposition~\ref{main-prop:singleleft}]
Using $\Atil=A+D\ell^+\ell^{+\top}/\kappa$ and Lemma~\ref{lem:bilinear},
\begin{equation}
\ell^{+\top}\Atil=\ell^{+\top}A+\frac{(\ell^{+\top}D\ell^+)}{\kappa}\,\ell^{+\top}
=\lambda^+\ell^{+\top}-2\lambda^+\ell^{+\top}=-\lambda^+\ell^{+\top}.
\end{equation}
For $\widetilde\ell_i^-=\ell_i^--c_i\ell^+$ with $c_i\coloneq\ell_i^{-\top}\Sigma\ell^+/\kappa$,
\begin{equation}
\widetilde\ell_i^{-\top}\Atil
=\lambda_i^-\ell_i^{-\top}-c_i\lambda^+\ell^{+\top}
+\frac{\ell_i^{-\top}D\ell^+}{\kappa}\ell^{+\top}
-c_i\frac{\ell^{+\top}D\ell^+}{\kappa}\ell^{+\top}.
\end{equation}
Substituting both parts of Lemma~\ref{lem:bilinear} and $c_i=\ell_i^{-\top}\Sigma\ell^+/\kappa$, the coefficient of $\ell^{+\top}$ collapses:
\begin{equation}
-c_i\lambda^+-(\lambda_i^-+\lambda^+)c_i+2c_i\lambda^+=-\lambda_i^-c_i,
\end{equation}
leaving $\widetilde\ell_i^{-\top}\Atil=\lambda_i^-\ell_i^{-\top}-\lambda_i^-c_i\ell^{+\top}
=\lambda_i^-\widetilde\ell_i^{-\top}$.
\end{proof}

\begin{proof}[Proof of Proposition~\ref{main-prop:committor}]
Under the condition of this proposition, $-\lambda^+$ is the eigenvalue of $\Atil$ closest to zero with eigenvector pair $\{\widetilde r^+,\widetilde\ell^+\}$, so for large $\tau$, $e^{\Atil\tau}\approx e^{-\lambda^+\tau}\widetilde r^+\widetilde\ell^{+\top}$, with the other terms suppressed by $e^{(\lambda_i^-+\lambda^+)\tau}\to0$. Then $K=\Stil^{1/2}e^{\Atil^{\!\top}\tau}\Stil^{-1/2}$ is asymptotically rank one, $K\approx e^{-\lambda^+\tau}(\Stil^{1/2}\widetilde\ell^+)(\Stil^{-1/2}\widetilde r^+)^{\top}$, whose SVD has leading left singular vector $\propto\Stil^{1/2}\widetilde\ell^+$; hence $v_1=\Stil^{-1/2}(\Stil^{1/2}\widetilde\ell^+)\propto\widetilde\ell^+\propto\ell^+$.
\end{proof}

\begin{proof}[Proof of Proposition~\ref{main-prop:singlesigma}]
Since $\Atil$ is Hurwitz (Proposition~\ref{main-prop:reflect}), the Lyapunov equation $\Atil\Stil+\Stil\Atil^{\!\top}=-D$ has a \emph{unique} solution; it therefore suffices to verify that Eq.~\eqref{main-eq:singlesigma} satisfies it. Write $u\coloneq\Sigma\ell^+$ and $\kappa\coloneq\ell^{+\top}\Sigma\ell^+$, so that $P_+=\ell^+\ell^{+\top}/\kappa$, $\Atil=A+D\ell^+\ell^{+\top}/\kappa$, and $\Stil=\Sigma-ur^{+\top}-r^+u^{\top}$.

Now expand $\Atil\Stil+\Stil\Atil^{\!\top}$ into four pieces coming from $\Atil=A+DP_+$ and $\Stil=\Sigma-ur^{+\top}-r^+u^{\top}$:
\begin{equation}
\Atil\Stil+\Stil\Atil^{\!\top}
=\underbrace{\big(A\Sigma+\Sigma A^{\!\top}\big)}_{(\mathrm I)}
+\underbrace{\big(DP_+\Sigma+\Sigma P_+D\big)}_{(\mathrm{II})}
-\underbrace{\Atil\big(ur^{+\top}+r^+u^{\top}\big)}_{(\mathrm{III})}
-\underbrace{\big(ur^{+\top}+r^+u^{\top}\big)\Atil^{\!\top}}_{(\mathrm{IV})}.
\end{equation}
Term (I) is $-D$, the target. We show (II)$-$(III)$-$(IV)$=0$.

For (II), $P_+\Sigma=\ell^+\ell^{+\top}\Sigma/\kappa=\ell^+u^{\top}/\kappa$, so $DP_+\Sigma=D\ell^+u^{\top}/\kappa$ and
\begin{equation}
(\mathrm{II})=\frac1\kappa\big(D\ell^+u^{\top}+u\ell^{+\top}D\big).
\end{equation}

For (III), use $\Atil u=\Atil\Sigma\ell^+=-\lambda^+\Sigma\ell^+=-\lambda^+u$ (Proposition~\ref{main-prop:reflect}) and $\Atil r^+=Ar^++DP_+ r^+$. Because $\ell^{+\top}r^+=1$, $\Atil r^+ = A r^+ + D\ell^+(\ell^{+\top}r^+)/\kappa =\lambda^+r^++D\ell^+/\kappa$. Hence
\begin{equation}
(\mathrm{III})=\Atil u\,r^{+\top}+\Atil r^+u^{\top}
=-\lambda^+ur^{+\top}+\Big(\lambda^+r^++\tfrac1\kappa D\ell^+\Big)u^{\top}
=\tfrac1\kappa D\ell^+u^{\top}+\lambda^+(r^+u^\top-ur^{+\top}).
\end{equation}
Term (IV) is the transpose of (III):
\begin{equation}
(\mathrm{IV})=\tfrac1\kappa u\ell^{+\top}D-\lambda^+(r^+u^\top-ur^{+\top}).
\end{equation}

Thus, $(\mathrm{II})-(\mathrm{III})-(\mathrm{IV})=0$, and $\Atil\Stil+\Stil\Atil^{\!\top}=-D$. By uniqueness, Eq.~\eqref{main-eq:singlesigma} is the stabilized covariance.
\end{proof}
\emph{Remark.} An equivalent, constructive proof is the spectral sum $\Stil=\int_0^\infty e^{\Atil t}De^{\Atil^{\!\top}t}\,dt$, evaluated block-by-block in the $\{\widetilde\ell^+,\widetilde\ell_i^-\}$ basis using Proposition~\ref{main-prop:singleleft}; both routes agree.


\section{Linear committors and the unstable left eigenvector}
\label{app:linear_committor}

We clarify the relation between the unstable left eigenvector and committor functions for a linear saddle. This relation is exact when the basin boundaries are compatible with a single linear coordinate; for general basin geometries, the committor need not have this one-dimensional form.

Consider the OU process centered at the fixed point,
\begin{equation}
    \dot{x} = A x + \xi,
    \qquad
    \langle \xi(t)\xi(t')^\top\rangle
    = D\delta(t-t'),
\end{equation}
with the backward generator
\begin{equation}
    L^\dagger
    =
    (Ax)\cdot\nabla
    +
    \frac{1}{2}D^{\alpha\beta}
    \partial_\alpha\partial_\beta.
\end{equation}
Let $\theta\neq 0$ and consider a nonconstant function depending on the single linear coordinate
\begin{equation}
    z=\theta^\top x,
    \qquad
    q(x)=Q(z).
\end{equation}
Its derivatives are
\begin{equation}
    \nabla q=\theta Q'(z),
    \qquad
    \partial_\alpha\partial_\beta q
    =
    \theta_\alpha\theta_\beta Q''(z),
\end{equation}
and hence the harmonic equation $L^\dagger q=0$ becomes
\begin{equation}
    \langle x,A^\top\theta\rangle Q'(z)
    +
    \frac{1}{2}
    \langle\theta,D\theta\rangle Q''(z)
    =
    0.
    \label{eq:linear_coordinate_generator}
\end{equation}

Assuming that $Q'$ is not identically zero, Eq.~\eqref{eq:linear_coordinate_generator} can be written as
\begin{equation}
    d\log Q'(z)
    =
    -\frac{2}{\langle\theta,D\theta\rangle}
    \langle x,A^\top\theta\rangle\,\theta^\top dx.
    \label{eq:committor_path_integral}
\end{equation}
The left-hand side is an ordinary one-dimensional differential, whereas the right-hand side is the line integral of a vector field in the full state space. For $Q'$ to be a well-defined function of $z=\theta^\top x$, this line integral must be independent of the integration path. The corresponding vector field must therefore have vanishing curl:
\begin{equation}
\begin{aligned}
    \partial_i
    \left[
        \langle x,A^\top\theta\rangle\theta_j
    \right]
    &=
    \partial_j
    \left[
        \langle x,A^\top\theta\rangle\theta_i
    \right]
    \\
    \Longleftrightarrow\qquad
    (A^\top\theta)_i\theta_j
    &=
    (A^\top\theta)_j\theta_i
    \\
    \Longleftrightarrow\qquad
    A^\top\theta\,\theta^\top
    &=
    \theta\theta^\top A.
    \label{eq:vanishing_curl}
\end{aligned}
\end{equation}
Because $\theta\neq 0$, Eq.~\eqref{eq:vanishing_curl} implies
\begin{equation}
    A^\top\theta
    =
    \frac{\theta^\top A\theta}
         {\theta^\top\theta}\,\theta
    =\lambda\theta.
    \label{eq:committor_left_eigenvector}
\end{equation}
Thus, a nonconstant harmonic function of the form $Q(\theta^\top x)$ can exist only if $\theta$ is a left eigenvector of $A$. Conversely, when $A^\top\theta=\lambda\theta$, the harmonic equation reduces to
\begin{equation}
    \lambda zQ'(z)
    +
    \frac{1}{2}d_\theta Q''(z)
    =
    0,
    \qquad
    d_\theta\coloneq\theta^\top D\theta,
    \label{eq:committor_ode}
\end{equation}
and therefore
\begin{equation}
    Q'(z)
    =
    C
    \exp\left(
        -\frac{\lambda z^2}{d_\theta}
    \right).
    \label{eq:committor_derivative}
\end{equation}

To obtain a committor, consider two absorbing hyperplanes
\begin{equation}
    \Gamma_A
    =
    \{x:\theta^\top x=z_A\},
    \qquad
    \Gamma_B
    =
    \{x:\theta^\top x=z_B\},
    \qquad
    z_A<z_B,
\end{equation}
and define
\begin{equation}
    q(x)
    =
    \mathbb{P}_x(\tau_B<\tau_A),
\end{equation}
where $\tau_A$ and $\tau_B$ are the respective first hitting times. The boundary conditions
\begin{equation}
    Q(z_A)=0,
    \qquad
    Q(z_B)=1
\end{equation}
give
\begin{equation}
    q(x)
    =
    \frac{
        \displaystyle
        \int_{z_A}^{\theta^\top x}
        \exp\left(
            -\frac{\lambda s^2}{d_\theta}
        \right)\,ds
    }{
        \displaystyle
        \int_{z_A}^{z_B}
        \exp\left(
            -\frac{\lambda s^2}{d_\theta}
        \right)\,ds
    }.
    \label{eq:linear_committor}
\end{equation}
Thus, every isocommittor is a hyperplane $\theta^\top x=\mathrm{const.}$ normal to the left eigenvector $\theta$.

For a saddle with a single unstable mode, take
$\theta=\ell^+$ and $\lambda=\lambda^+>0$. Writing
\begin{equation}
    \alpha
    \coloneq
    \sqrt{
        \frac{\lambda^+}
             {\ell^{+\top}D\ell^+}
    },
\end{equation}
the committor between the finite absorbing hyperplanes $\ell^{+\top}x=z_A$ and $\ell^{+\top}x=z_B$ is
\begin{equation}
    q(x)
    =
    \frac{
        \operatorname{erf}\!\left(
            \alpha\ell^{+\top}x
        \right)
        -
        \operatorname{erf}\!\left(
            \alpha z_A
        \right)
    }{
        \operatorname{erf}\!\left(
            \alpha z_B
        \right)
        -
        \operatorname{erf}\!\left(
            \alpha z_A
        \right)
    }.
    \label{eq:finite_error_function_committor}
\end{equation}
Thus, the committor has an affinely normalized error-function profile along the unstable coordinate, and its isocommittors are hyperplanes normal to $\ell^+$. In the far-field limit $z_A\to-\infty$ and $z_B\to+\infty$, this reduces to
\begin{equation}
    q(x)
    =
    \frac{1}{2}
    \left[
        1+
        \operatorname{erf}\!\left(
            \alpha\ell^{+\top}x
        \right)
    \right],
\end{equation}
which is the error-function committor of the harmonic-saddle approximation~\cite{main-berezhkovskii2005one}.

The restriction on the boundary geometry is essential. For general reactant and product basins, the boundary conditions can introduce dependence on the transverse coordinates even when the dynamics is linear.
